\documentclass{article}
\usepackage{amsmath, amssymb, xcolor, amsthm}
\usepackage{bbm}
\usepackage{graphicx,pdfsync}
\usepackage{amsxtra,amssymb,amsmath,enumerate,graphicx,bbm,hyperref,amsthm,color}
\usepackage[active]{srcltx}
\usepackage{t1enc}
 
\usepackage[T1]{fontenc}

\begin{document}
\pagestyle{plain}
\bibliographystyle{plain}
\newtheorem{theo}{Theorem}[section]
\newtheorem{lemme}[theo]{Lemma}
\newtheorem{example}[theo]{Example}
\newtheorem{cor}[theo]{Corollary}
\newtheorem{defi}[theo]{Definition}
\newtheorem{prop}[theo]{Proposition}
\newtheorem{problem}[theo]{Problem}
\newtheorem{remarque}[theo]{Remark}
\newtheorem{Exemple}[theo]{Exemple}
\newcommand{\beq}{\begin{eqnarray}}
\newcommand{\enq}{\end{eqnarray}}
\newcommand{\be}{\begin{eqnarray*}}
\newcommand{\en}{\end{eqnarray*}}
\newcommand{\ben}{\begin{eqnarray*}}
\newcommand{\enn}{\end{eqnarray*}}
\newcommand{\Td}{\mathbb T^d}
\newcommand{\Rd}{\mathbb R^n}
\newcommand{\R}{\mathbb R}
\newcommand{\N}{\mathbb N}
\newcommand{\Sn}{\mathbb S}
\newcommand{\Zd}{\mathbb Z^d}
\newcommand{\Linf}{L^{\infty}}
\newcommand{\dt}{\partial_t}
\newcommand{\Dt}{\frac{d}{dt}}
\newcommand{\Dtt}{\frac{d^2}{dt^2}}
\newcommand{\demi}{\frac{1}{2}}
\newcommand{\vf}{\varphi}
\newcommand{\epu}{_{\varepsilon}}
\newcommand{\ep}{^{\varepsilon}}
\newcommand{\bfi}{{\mathbf \Phi}}
\newcommand{\bpsi}{{\mathbf \Psi}}
\newcommand{\bx}{{\mathbf x}}
\newcommand{\dis}{\displaystyle}
\newcommand{\ds}{\partial_s}
\newcommand{\dss}{\partial^{2}_{s}}
\newcommand{\dx}{\partial_x}
\newcommand{\dxx}{\partial^{2}_{x}}
\newcommand{\dy}{\partial_y}
\newcommand{\dyy}{\partial^{2}_{y}}
\newcommand{\dtt}{\partial^{2}_{t}}
\newcommand {\g}{\`}
\newcommand{\E}{\mathbb E}
\newcommand{\bQ}{\mathbb Q}
\newcommand{\1}{{\mathbf I}}
\newcommand{\bF}{\mathbb F}
\newcommand{\F}{\cal F}
\newcommand{\bP}{\mathbb P}
\let\cal=\mathcal
\newcommand{\lb}{\langle}
\newcommand{\rb}{\rangle}
\newcommand{\bv}{\bar{v}}
\newcommand{\uv}{\underline{v}}
\newcommand{\cS}{{\cal S}}
\newcommand{\cSf}{\cS_\Phi}
\newcommand{\bw}{\bar{w}}
\newcommand{\uw}{\underline{w}}
\newcommand{\bsig}{\bar{\sigma}}
\newcommand{\usig}{\underline{\sigma}}
\newcommand{\ms}{\mathfrak{\sigma}}
\newcommand{\mm}{\mathfrak{\mu}}
\newcommand{\Ac}{\mathcal{A}} 
\newcommand{\Fc}{\mathcal{F}}  
\newcommand{\Dc}{\mathcal{D}} 
\newcommand{\Gc}{\mathcal{G}}  
\newcommand{\vr}{{\rm v}}
 \newcommand{\vrt}{\tilde{\rm v}} 
\newcommand{\vrb}{{\rm \underline v}} 
\newcommand{\eps}{\varepsilon}  
\newcommand{\Pas}{\bP-{\rm a.s.}}  
\newcommand{\x}{\times}  
\def\red#1{{\color{red} #1}}
\def\blue#1{{\color{blue} #1}}
\def\bru#1{{\color{green} #1}}
\def \endp{\hbox{ }\hfill$\Box$}
\def \vp{\varphi}
\parindent=0mm
\def\vs#1{\vspace{#1mm}}
\def\sr{{\rm s}}
 \def\bS{{\rm S}}
  \def\smf{{\mathfrak s}}

\author{Bruno {\sc Bouchard}\footnote{Universit\'e Paris-Dauphine, PSL  University, CNRS, CEREMADE, 75016 Paris, France, bouchard@ceremade.dauphine.fr. ANR Grant CAESARS (ANR-15-CE05-0024), Initiative de Recherche ``M\'ethodes non-lin\'eaires pour la gestion des risques financiers'' sponsored by AXA Research Fund.}, Gr\'egoire {\sc Loeper}\footnote{Monash University, School of Mathematical Sciences \& Centre for Quantitative Finance and Investment Strategies (CQFIS),  gregoire.loeper@monash.edu. CQFIS has been supported by BNP Paribas.}, \\ Halil Mete {\sc Soner}\footnote{ETH Zurich, mete.soner@math.ethz.ch. Partially supported by the ETH Foundation, Swiss Finance Institute and the Swiss National Foundation through SNF 200020-172815.}, Chao {\sc Zhou}\footnote{Department of Mathematics, National University of Singapore, Singapore, matzc@nus.edu.sg. Research supported by  Singapore MOE AcRF Grants R-146-000-219-112 and R-146-000-255-114.}} 
\title{Second order stochastic target problems with generalized market impact}
\date{\today}
 
 \maketitle
 
\begin{abstract}
We extend the study of \cite{BoLoZo2,LoMi1} to  stochastic target problems with general market impacts. Namely, we consider a general abstract model which can be associated to a fully nonlinear parabolic equation. Unlike \cite{BoLoZo2,LoMi1}, the equation is not concave and the regularization/verification approach of \cite{BoLoZo2} can not be applied.  {We also relax the gamma constraint of \cite{BoLoZo2}}. In place, we need to generalize the a priori estimates   of \cite{LoMi1} and exhibit smooth solutions from the classical parabolic equations theory. Up to an additional approximating argument, this allows us to show that the super-hedging price solves the parabolic equation and that a perfect hedging strategy can be constructed when the coefficients are smooth enough. This representation leads to a general  dual formulation. We finally provide an asymptotic expansion around a model without impact.  
\end{abstract}

\section{Introduction}
 Inspired by \cite{AbergelLoeper,LoMi1},  {the authors in} \cite{BoLoZo1,BoLoZo2} considered a financial market with permanent price impact (and possibly a resilience effect), in which the impact function behaves as a linear function (around the origin) in the number of purchased stocks.  This class of models is dedicated to the pricing and hedging of derivatives in situations  where the notional of the product hedged is such that the delta-hedging is non-negligible compared to the average daily volume traded on the underlying asset. As opposed to \cite{BoLoZo1}, the options considered in \cite{BoLoZo2,LoMi1} are covered, meaning that  the buyer of the option delivers, at the inception, the required initial delta position, and accepts a mix of stocks (at their current market price) and cash as payment for the final claim. This is a common practice which eliminates the cost incurred by the initial and final hedge. In \cite{LoMi1}, the author considers a Black-Scholes type model, while the model of \cite{BoLoZo2} is a local volatility one.

Motivated by these works, we consider in this paper a general abstract model of market impact in which the dynamics of the stocks $X$, the  {\sl wealth}\footnote{{More precisely: the value of the cash plus the number of stocks in the portofolio times the current value of the stocks.}} $V$ and the number of stocks $Y$ held in the portfolio follow dynamics  of the form 
\begin{align*}
X&=x+ \int_{t}^{\cdot} \mu(s,X_{s},\gamma_{s}{,b_{s}}) d{s} + \int_{t}^{\cdot} \sigma(s,X_{s},\gamma_{s}) dW_{s}\\
Y&=y+\int_{t}^{\cdot} b_{s}ds + \int_{t}^{\cdot} \gamma_{s} dX_{s}\\
V&=v+\int_{t}^{\cdot} F(s,X_{s},\gamma_{s})ds + \int_{t}^{\cdot} Y_{s} dX_{s}
\end{align*}
where $(y,b,\gamma)$ are the controls, and we consider the general super-hedging problem: 
 \beq
\vr(t,x):=\inf\{ v=c+ yx~:~(c,y)\in \R^{2} \mbox{ s.t. }    \Gc(t,x,v,y)\ne \emptyset\},\nonumber
\enq
in which 
\beq 
\Gc(t,x,v,y)= \Big\{(b,\gamma) :  
 V^{t,x,v,\phi}_{T}\ge g(X^{t,x,\phi}_{T})\mbox{ for   } \phi:=(y,b,\gamma)\Big\}\nonumber,
\enq
and $g$ is the payoff function associated to a European claim.   

One can easily be convinced, by using formal computations based on the geometric dynamic programming principle of \cite{SonTouzDyn}, see also the discussion just after Remark \ref{rem: ex bar F}, that $\vr$ should be a super-solution of  the fully nonlinear parabolic equation 
$$
0\le -\partial_{t}\vr-\bar F(\cdot,\partial^{2}_{x} \vr)\;\;  \mbox{ and } \;\; {(|F|+|\sigma|)}(\cdot,\partial^{2}_{x} \vr)<\infty.
$$
in which  
 $$
\bar F(t,x,z):= \frac12 \sigma(t,x,z)^{2} z - F(t,x,z).
$$
The right-hand side constraint in the previous inequalities is of importance. Indeed ${(F,\sigma)}(t,x,\cdot)$ can typically be singular and only finite on an interval of the form $(-\infty,\bar \gamma(t,x))$, as it is the case in  
 \cite{BoLoZo2}. Under this last assumption, one can actually expect that $\vr$ is a viscosity solution of  
 \begin{align}\label{eq: pde vr intro}
 \min\{-\partial_{t}\vr-\bar F(\cdot,\partial^{2}_{x} \vr)\;,\; \bar \gamma - \partial^{2}_{x}\vr\}=0\mbox{ on } [0,T)\x \R,
  \end{align} 
  with $T$-terminal condition given by the smallest function $\hat g\ge g$ such that $\dxx \hat g\le \bar \gamma(T,\cdot)$.
 
 In  \cite{BoLoZo2}, the authors impose a strong  (uniform) constraint on the controls of the form $\gamma\le \tilde  \gamma(\cdot,X^{t,x,\phi})$ with $\tilde \gamma$ such that $F(\cdot,\tilde \gamma)\le C$ for some $C>0$, and obtain that $\vr$ is actually the unique viscosity solution of \eqref{eq: pde vr intro} with $\tilde \gamma$ in place of $\bar \gamma$, and terminal condition $\hat g$ (defined with $\tilde \gamma$ as well).  Their proof of the super-solution property mimicks   arguments  of \cite{CheriSonTouz},  and we can follow this approach. As for the sub-solution property,  they could not prove  the appropriate dynamic programming principle, and  the standard direct arguments could not be used. Instead,  they employed a regularization argument for viscosity solutions, inspired by \cite{krylov2000rate}, together with a verification procedure. In \cite{BoLoZo2}, the authors critically use the fact that $\bar F$ is convex. 
 
Our setting here is different. First, as in \cite{LoMi1}, we do not impose a uniform constraint on our strategies. Our controls can take values arbitrarily close to the singularity $\bar \gamma(\cdot,X^{t,x,\phi})$ and the equation \eqref{eq: pde vr intro} is possibly degenerate. Even for $\bar F$ defined as in \cite{BoLoZo2} our setting is more general in a sense.  Second, $\bar F$ is not assumed to be convex. 

For these reasons, we can not reproduce the smoothing/verification argument of \cite{BoLoZo2} to deduce that $\vr$ is actually a subsolution.
\vs2

In this paper, we therefore proceed differently and generalise arguments used in \cite{LoMi1} in the context of a Black-Scholes type model. Namely, we directly use the theory of parabolic equations  to prove the existence of  smooth solutions to \eqref{eq: pde vr intro}
whenever $\hat g$ is smooth and satisfies a constraint of the form $\dxx \hat g\le \bar \gamma(T,\cdot)-\eps$, for some $\eps>0$.  Our analysis heavily relies on new {\it a priori} estimates, see Proposition \ref{prop: u in DC bar gamma eps} below, thanks to which one can appeal to the continuity method in a rather classical way, see the proof of Theorem \ref{theo: u C4}. We then let $\eps$ go to $0$ to conclude that $\vr$ indeed solves \eqref{eq: pde vr intro} in the viscosity solution sense, see Theorem \ref{theo: vr=u for Phi=hat g 0} below.
\vs2

We also discuss two important issues that were not considered in \cite{BoLoZo2} but already studied in \cite{LoMi1} in a Black-Scholes type model: 

- The first one concerns the asymptotic expansion of the price around a model without market impact. As in  \cite{LoMi1}, we show that a first order expansion can be established, see Proposition \ref{prop: asymptotique} below. But, we also prove that one can deduce from it a strategy that matches the terminal face-lifted payoff $\hat g$ at any prescribed level of precision in ${\mathbb L}^{\infty}$-norm, see Proposition \ref{prop: approximate hedging}.  

- The second one concerns the existence of a dual formulation. It can be established when $\bar F$ is convex in its last argument, see Theorem \ref{thm: dual formulation}. Applied to  the model discussed in \cite{BoLoZo2},  see Example \ref{exemple: BoLoZo} below, it takes the form 
\begin{align*}
 \vr(t,x)&=\sup_{\smf } \E\left[\hat g(X^{t,x,\smf}_{T})-\int_{t}^{T} \frac12 \frac{(\smf_{s}-\sigma_\circ(t,X^{t,x,\smf}_{s}))^{2}}{f(X^{t,x,\smf}_{s})} ds\right]\\
 &=\sup_{\smf } \E\left[ g(X^{t,x,\smf}_{T})-\int_{t}^{T} \frac12 \frac{(\smf_{s}-\sigma_\circ(t,X^{t,x,\smf}_{s}))^{2}}{f(X^{t,x,\smf}_{s})} ds\right]
\end{align*}
 in which 
  $
 X^{t,x,\smf}=x+\int_{t}^{\cdot} \smf_{s}dW_{s}
 $, $\sigma_{\circ}$ is the volatility surface in a the market without impact and $f>0$ is the impact function, the limit case $f\equiv 0$ corresponding to the absence of  impact.  It can be interpreted as the formulation of the super-hedging price with volatility uncertainty. The difference being that the formula is penalized by the squared distance of the realized volatility  term $\smf$ to the original local volatility  $\sigma_\circ(\cdot,X^{t,x,\smf})$ associated to the model, weighted by the inverse of the impact function $f(X^{t,x,\smf})$. It can also be seen as a martingale  optimal transport problem, see  \cite[Section 4.1]{LoMi1}
 for details. \vs2
 
 To conclude, let us refer to \cite{becherer2016optimalbis,becherer2016optimal,becherer2017stability,CetinJarrowProtter,CheriSonTouz,Frey,Liu,Schon,Sircar,SonTouzDyn}, and the references therein. Also for related works, see \cite{BoLoZo2} for a discussion.
 
 \vs2
 
The rest of this paper is organized as follows. The general abstract market model is described in Section \ref{sec: model} and the characterization of $\vr$ as a solution of a parabolic equation is proved in Section \ref{sec: pde cara}. The asymptotic expansion and the dual formulation are provided and discussed in Sections \ref{sec: asympto} and \ref{sec: dual}.
\\

\noindent{\bf General notations.}
Throughout this paper,  $\Omega$ is the canonical space of continuous functions on $\R_{+}$ starting at $0$, $\bP$ is the Wiener measure, $W$ is the canonical process, and $\F=(\Fc_{t})_{t\ge 0}$ is {the augmentation of its} raw filtration $\F^{\circ}=(\Fc^{\circ}_{t})_{t\ge 0}$.  All random variables are defined on $(\Omega,\Fc_{\infty},\bP)$. We denote by $|x|$   the Euclidean norm of $x\in \R^{n}$, the integer  $n\ge 1$ is given by the context. {Unless otherwise specified, inequalities involving random variables are taken in the $\Pas$ sense. We use the convention $x/0={\rm sign}(x)\times \infty$ with ${\rm sign}(0)= +$. We denote by $\partial^{n}_{x}\vp$ the $n$th-order derivative of a function $\vp$ with respect to its $x$-component, whenever it is well-defined. } For $E,F,G$, three subsets of $\R$, We denote by $C^{h,k}_b(E\times F)$ the set of continuous functions on $E\times F$ which have bounded partial derivatives of order from 1 to $h$ with respect to the first variable and from 1 to $k$ to the second variable. We denote by $C^{{h,k,l}}(E\times F\times G)$ the set of continuous functions on $E\times F\times G$ which have partial derivatives of order from 1 to $h$ with respect to the first variable, from 1 to $k$ to the second variable and from 1 to $l$ to the third variable. We denote by $C^{h}_b(E\times F)$ the set of continuous functions on $E\times F$ which have bounded partial derivatives of order   1 to $h$.   If in addition its $h$-th order derivatives are uniformly $\alpha$-H\"older, with $\alpha\in (0,1)$, we say that it belongs to $C^{h+\alpha}_{b}(E\times F)$. We omit the spaces $E,F,G$ if they are clearly given by the context. 

\section{Abstract market impact model}\label{sec: model}

We first describe our abstract  market with impact. It generalizes the model studied in \cite{BoLoZo1,BoLoZo2,LoMi1}. We use the representation of the hedging strategies described in \cite{BoLoZo2}, which is necessary to obtain the supersolution characterization of the super-hedging price of Proposition \ref{prop: visco sol vr eps} below. How to get to the market evolution (\ref{eq: def X}, \ref{eq: def Y}, \ref{eq: def V}) is explained  {briefly} in  {Example} \ref{exemple: BoLoZo}.
\\

More precisely, given $k\ge 1$, we denote by  ${\Ac^{\circ}_{k}}$ the collection of continuous and $\bF$-adapted processes $(b,\gamma)$ such that 
$$
\gamma=\gamma_{0}+\int_{0}^{\cdot}\beta_{s} ds +\int_{0}^{\cdot} \alpha_{s} dW_{s}
$$
where $(\alpha,\beta)$ is  continuous, $\bF$-adapted,  and $\zeta:=(b,\gamma,\alpha,\beta)$ is essentially bounded by $k$  and such that
\begin{align*}
\E\left[\sup\left\{|\zeta_{s'}-\zeta_{s}|,\; t\le s \le s'\le s+\delta\le T\right\}|\Fc_{t}^{\circ} \right]\le k\delta
\end{align*}
for all $0\le \delta\le 1$ and $t\in [0,T-\delta]$. We then define 
\begin{align*} 
& {\Ac^{\circ}:=\cup_{k}\Ac^{\circ}_{k}}.&
\end{align*}

Let $F: [0,T]\x \R^{2}\mapsto \R\cup\{\infty\}$ be a continuous map and let 
$$
\Dc:=\{F<\infty\}
$$
be its domain.  We assume that there exists a  map $(t,x) \to \bar\gamma(t,x)\in \R\cup\{+\infty\}$  such that 
\beq\label{eq:defbargamma}
\Dc=\{{(t,x,z)\in [0,T]\x \R \x \R:  \ z\in (-\infty, \bar\gamma(t,x))}\},
\enq
and  that 
\begin{equation}\label{eq: hyp bar gamma unif cont}
\mbox{$\bar\gamma$ is either uniformly continuous, or identically equal to $+\infty$.}
\end{equation}


We now let $\mu:\Dc\x\R\to \R$ and $\sigma: \Dc\to \R$ be two continuous maps such that, for all $\eps>0$,

\begin{align}\label{eq: hyp sigma}
&&\mbox{$\mu$ is Lipschitz, with linear growth in its second variable, on  $\Dc_{\eps,\eps^{-1}}\x \R$},\\ 
&&\mbox{ $\sigma$  is Lipschitz, with linear growth in its second variable, on $\Dc_{\eps,\eps^{-1}}$,} \nonumber
\end{align}

where  
\begin{eqnarray}\label{defDeps}
&&\Dc_{\eps}:=\{(t,x,z)\in [0,T]\x \R^{2}: F(t,x,z)\le \eps^{-1}\},\\ 
&&\Dc_{\eps,k}:=\Dc_{\eps}\cap ([0,T]\x\R\times [-k,k]) \mbox{ for } k>0.\nonumber
\end{eqnarray}

Then, given $(t,x,v)\in [0,T]\x \R\x \R$ and  $\phi=(y,b,\gamma)\in \R\x \Ac^{\circ}$, we define $(X^{t,x,\phi},Y^{t,x,\phi},$ $V^{t,x,v,\phi})$ as the solution  on $[t,T]$ of 
\begin{align}
X&=x+ \int_{t}^{\cdot} \mu(s,X_{s},\gamma_{s}{,b_{s}}) d{s} + \int_{t}^{\cdot} \sigma(s,X_{s},\gamma_{s}) dW_{s}\label{eq: def X}\\
Y&=y+\int_{t}^{\cdot} b_{s}ds + \int_{t}^{\cdot} \gamma_{s} dX_{s}\label{eq: def Y}\\
V&=v+\int_{t}^{\cdot} F(s,X_{s},\gamma_{s})ds + \int_{t}^{\cdot} Y_{s} dX_{s}\label{eq: def V}
\end{align}
satisfying $(X_{t},Y_{t},V_{t})=(x,y,v)$, whenever $(\cdot,X,\gamma)$ takes values in $\Dc_{\eps,k}$ 
on $[0,T]$, for some $\eps,k>0$. 
If this is the case, we say that $\phi$ belongs to $\Ac_{k}^{\eps}$. For ease of notations, we set $\Ac:=\cup_{\eps,k>0} \Ac_{k}^{\eps}$. 

For a payoff function  $g:\R\to \R$ the super-hedging price of the covered European claim associated to $g$    is then defined as 
\beq
\label{defv}\vr(t,x):=\inf\{ v=c+ yx~:~(c,y)\in \R^{2} \mbox{ s.t. }    \Gc(t,x,v,y)\ne \emptyset\},
\enq
in which 
\beq 
\Gc(t,x,v,y)= \Big\{\phi=(y,b,\gamma)\in \Ac:  
 V^{t,x,v,\phi}_{T}\ge g(X^{t,x,\phi}_{T})\Big\}\nonumber
\enq
whenever this set is non-empty.  Note that 
 \begin{align}\label{eq: def v veps vepsk}
 \vr(t,x)=\inf_{\eps>0}\vr^{\eps}(t,x)\; \mbox{ where }\; \vr^{\eps}(t,x):=\inf_{k>0}\vr^{\eps}_{k}(t,x)
 \end{align}
 in which $\vr^{\eps}_{k}$ is defined as $\vr$ but in terms of $\Ac_{k}^{\eps}$. 
 \vs2
 
 In the following, we assume as in \cite{BoLoZo2} that 
 \begin{align}\label{eq: hyp g}
 \mbox{$g$ is lower-semicontinuous, bounded from below, and  $g^{+}$ has linear growth.}
 \end{align}

\begin{example}[Example of derivation of the evolution equations]\label{exemple: BoLoZo}

{We close this section with an example of formal derivation of the above abstract dynamics. }
In the spirit of \cite{AbergelLoeper,LoMi1}, {let us consider a   linear market impact model in which an (infinitesimal) order to buy $dY_{t}$  stocks at $t$ leads to a permanent price move of $ f(t,X_{t},\gamma_{t})dY_{t}$, and to an average execution price of $X_{t}+ f(t,X_{t},\gamma_{t})dY_{t}  +  \bar f(t,X_{t},\gamma_{t})dY_{t}$. Then, following the computations done in \cite{AbergelLoeper,LoMi1}, see also the rigorous proof in \cite{BoLoZo1} for details\footnote{ The continuous time version is obtained by considering the limit dynamics of a discrete time trading model, as the speed of trading goes to infinity.}, the {\sl portfolio value} $V$ corresponding to the holding in cash plus the number of stocks in the portfolio evaluated at their current price $X$ is given by\footnote{Obviously, this is only a theoretical value, the liquidation value of the portfolio being different.}
}
\ben
V=v + \int_t^\cdot Y_s dX_s-\int_t^\cdot {\bar f}({s,X_{s}},\gamma_{s}) {d\lb Y,Y \rb_s}.
\enn

The contribution ${\bar f}(s,X_{s},\gamma_{s}) d\lb Y,Y \rb_s$ is the spread between the execution price of the trade and the final price after market impact. It can be either positive or negative. The fact that $f$ and $\bar f$ can depend on $\gamma$ {is} discussed in \cite{LoMi1}. 

{Let us now assume that $X$ evolves according to $dX_t =  \sigma_\circ(t,X_{t}) dW_t +  \mu_\circ(t,X_{t}) dt$ in the absence of trade. Then, arguing again as in \cite{BoLoZo1}, }
we obtain  the modified dynamics
\ben
{dX_t =  \sigma_\circ(t,X_{t})  dW_t +  \mu_\circ(t,X_{t})  dt + f(t,X_t, \gamma_t) dY_t+ f'(t,X_t, \gamma_t) \gamma_{t} \sigma_\circ(t,X_{t})^{2}dt.}
\enn

Combining this with (\ref{eq: def Y}), {and formally solving in $dX$,}  we obtain that
\ben
\sigma(t,X_{t},\gamma_{t}) = \frac{ {\sigma_\circ}(t,X_{t})}{1-f(t,X_{t},\gamma_{t})\gamma_{t}},
\enn
{so that the dynamics of  $V$ can be written as}
\ben
V=v + \int_t^\cdot Y_s dX_s-\int_t^\cdot \bar f(s,X_{s},\gamma_{s}) \left(\frac{ {\sigma_\circ}(s,X_{s})\gamma_{s}}{1-f(s,X_{s},\gamma_{s})\gamma_{s}}\right)^2ds.
\enn
Note that, as observed in \cite{BoLoZo1}, the drift $  \mu_\circ$ is also affected by the market impact, but {that this} does not affect the pricing equation. It is therefore not taken into account in our abstract model.

 The model studied in \cite{BoLoZo1,BoLoZo2} corresponds to $f=f(x)$ (no dependency in $t,\gamma$) and  $ \bar f=-f/2$. In this particular case, the functions $\sigma$ and $F$ are given by
\begin{eqnarray*}
&\sigma(t,x,{z})=\frac{ \sigma_\circ(t,x)}{1-f(x)z},\; \bar \gamma = 1/f&\\
& F(t,x,z)=\frac12 \left(\frac{ \sigma_\circ(t,x){z}}{1-f(x)z}\right)^{2}f(x)\1_{\{f(x)z<1\}}{+\infty\1_{\{f(x)z\ge 1\}}}.&
\end{eqnarray*}
 
\end{example}


\section{PDE characterization}\label{sec: pde cara}

The parabolic equation associated to $\vr$ can be formally derived as follows. Assume that $\vr$ is smooth and that a perfect hedging strategy $\phi=(y,b,\gamma)$ can be found when starting at $t$ from $v=\vr(t,x)$ if the stock price is $x$ at $t$. Then, we expect to have $V^{t,x,v,\phi}=\vr(\cdot,X^{t,x,\phi})$ which, by It\^{o}'s lemma combined with \eqref{eq: def X}-\eqref{eq: def V}, implies that 
\begin{align*}
&F(s,X^{t,x,\phi}_{s},\gamma_{s})ds+Y^{t,x,\phi}_{s}dX^{t,x,\phi}_{s}\\
&=(\partial_{t}\vr+\frac12\sigma^{2}(\cdot,\gamma_{s})\partial^{2}_{x} \vr)(s,X^{t,x,\phi}_{s}) ds + \partial_{x} \vr(s,X^{t,x,\phi}_{s})dX^{t,x,\phi}_{s}
\end{align*}
for $s\in [t,T]$. By identifying the different terms, we obtain 
$$
F(s,X^{t,x,\phi}_{s},\gamma_{s})=(\partial_{t}\vr+\frac12\sigma^{2}(\cdot,\gamma_{s})\partial^{2}_{x} \vr)(s,X^{t,x,\phi}_{s})\mbox{ and }Y^{t,x,\phi}_{s}=\partial_{x} \vr(s,X^{t,x,\phi}_{s}).
$$
Another application of It\^{o}'s lemma to the second equation then leads to 
$$
\gamma_{s}=\partial^{2}_{x} \vr(s,X^{t,x,\phi}_{s}),
$$
recall \eqref{eq: def Y}.
The combination of the above reads 
$$
0=-(\partial_{t}\vr+\bar F(\cdot,\partial^{2}_{x} \vr))(s,X^{t,x,\phi}_{s}) \mbox{ and } { (|F|+|\sigma|)}(\cdot,\partial^{2}_{x} \vr)(s,X^{t,x,\phi}_{s})<\infty, 
$$
in which 
\begin{equation}\label{eq: def bar F}
\bar F(t,x,z):=\frac12 \sigma(t,x,z)^{2} z - F(t,x,z),\; \mbox{ for } (t,x,z)\in \Dc.
\end{equation}

\begin{remarque}\label{rem: ex bar F} The model discussed in {\cite{BoLoZo2}} corresponds to 
$$
 \bar F(t,x,z)=\frac12\frac{\sigma_\circ^{2}(t,x)z}{1-f(x)z}\1_{\{f(x)z<1\}} {+\infty\1_{\{f(x)z\ge 1\}}}.
 $$

\end{remarque}

As usual, perfect equality can not be ensured because of the gamma constraint induced by the above. We therefore only expect to have 
 $$
0\le -(\partial_{t}\vr+\bar F(\cdot,\partial^{2}_{x} \vr))(s,X^{t,x,\phi}_{s}) \mbox{ and }  {(|F|+|\sigma|)}(\cdot,\partial^{2}_{x} \vr)(s,X^{t,x,\phi}_{s})<\infty.
$$
Recalling  \eqref{eq:defbargamma}, this leads to the fact that $\vr$ should be a super-solution of the  parabolic equation 
 \begin{align}\label{eq: pde vr}
 \min\{-\partial_{t}\vp- \bar F(\cdot,\partial^{2}_{x}\vp)\;,\; \bar \gamma - \partial^{2}_{x}\vp\}=0\mbox{ on } [0,T)\x \R.
  \end{align}
  By minimality, it should indeed be a solution. Moreover, as usual, the gamma constraint $\partial^{2}_{x}\vp\le \bar \gamma$ needs to propagate up to the boundary, so that we can only expect that $\vr$ satisfies the terminal condition 
 \begin{align}\label{eq: pde vr T}
 \lim_{(t',x')\to (T,x)} \vp(t',x')=\hat g(x) \mbox{ for }  x\in \R, 
  \end{align}
  where $\hat g$ is the face-lift of $g$, i.e.
  $$
  \hat g=\inf\{h\in C^{2}(\R): h\ge g\mbox{ and } \partial^{2}_{x} h\le \bar \gamma(T,\cdot) \}. 
  $$
 See Remark \ref{rem: face-lift} below for ease of comparison with \cite{BoLoZo2}.
 
 \begin{remarque} When $\bar \gamma\equiv +\infty$, the above reads 
  \begin{align*}
 -\partial_{t}\vp- \bar F(\cdot,\partial^{2}_{x}\vp)=0\mbox{ on } [0,T)\x \R\;\mbox{ and }\;
 \lim_{(t',x')\to (T,x)} \vp(t',x')= g(x) \mbox{ for }  x\in \R.
  \end{align*} 
 \end{remarque}
 
 In order to prove that $\vr$ is actually a continuous viscosity solution of the above, we need some additional assumptions. First, we assume that $\bar F$ is smooth enough, 
\begin{align}
&  \bar F \in C^{{1}}(\Dc)\;\mbox{  and }  \; \bar F\in {C^{1,3,3}_b(\Dc_{\eps,\eps^{-1}})},\;\eps\in (0,\eps_{\circ}],\label{eq : regul bar F cas general}\\
& \ \bar F \mbox{ is uniformly continuous on } \Dc_{\eps}, \;\eps\in (0,\eps_{\circ}]\label{eq: bar F uniform cont sur Deps},
\end{align}
where $\eps_{\circ}>0$, 
and that 
\begin{align}
&F(\cdot,0)=0,\label{eq: F(,0)=0}.
\end{align}
For later use, note that the above implies  
\begin{align}
&\bar F(\cdot,0)=0 \label{eq: bar F(,0)=0}.
\end{align}
We also assume  that there exists $L_{\circ}, M>0$ such that, on $\Dc$ and for all $\eps\in (0,\eps_{\circ}]$,
\begin{align}
&|{\partial_{t}} \bar F/\bar F|\le L_{\circ}, \label{eq: bar Ft/ bar F bounded by Lo}  
\mbox{ and } 
|\partial^{2}_{x} \bar F(\cdot,z)| \leq M|z| \text{ for all } z\in (-\infty,0], 
\end{align}
that 
\begin{align}
&  \partial_{z} \bar F > 0\; \mbox{ on }\; \Dc_{\eps}\; \mbox{ and }\; \sup_{\{(t,x,z)\in \Dc_{\eps, \varepsilon^{-1}} 
\}}(|{\partial_{z}} \bar F|+ |1/{\partial_{z}} \bar F|)<\infty , \label{eq: bar F unif elliptic}\\
&\inf_{\Dc_{\eps, \varepsilon^{-1}}} \sigma>0. \label{eq: bound inf sigma}  \\
&\mbox{ $F$ is uniformly continuous on  $\Dc_{\eps}$},\label{eq: hyp F monotone unif cont}
\end{align}
and that, for all $\eps\in (0,\eps_{\circ}]$,  there exists a {continuous} map  $\bar\gamma_\eps$ such that
\begin{align}
& \Dc_{\eps}= \{(t,x,z)\in [0,T]\x \R^{2}: z\leq \bar\gamma_\eps(t,x)\} \; 
\end{align}

\begin{remarque} All these conditions are satisfied in the model of \cite{BoLoZo2}.
\end{remarque}
\begin{remarque} {As a corollary of (\ref{eq: F(,0)=0}) and (\ref{eq: hyp F monotone unif cont}), we have that
\beq
\sup_{   {\Dc_{\eps}}
}|  F|<\infty,  \label{eq: bound F sur Dc eps}
\enq}
\end{remarque}

Finally, we assume that   
\beq\label{eq: def hat g eps}
\hat{g}^{\eps}:=\inf\{h\geq g: h\in C^{2}(\R), {F}(T,\cdot, \dxx h) < \eps^{-1}\}\;
\enq  
satisfies 
\begin{align}\label{eq: hyp hat g eps}
\mbox{$\hat g^{\eps}$  is uniformly continuous, bounded from below and has linear growth} 
\end{align}
and that there exists $k_{\circ}\ge 1$ such that 
\begin{align}
&[\vr^{\eps}_{k}]^{+}   \mbox{ has linear growth, uniformly in $k\ge k_{\circ}$}, \label{eq: bound vr et vr eps k}
\end{align}
for all $0< \eps$ $\le$  $\eps_{\circ}  $, 
in which we use the convention $1/0=\infty$ and identify $\hat g$ with $\hat g^{0}$.
\vs2

Under the above conditions, we can state the main result of this section.
 
  \begin{theo}\label{theo: vr=u for Phi=hat g 0} The function $\vr$ is a continuous  viscosity solution of \eqref{eq: pde vr}
  such that  $\lim_{t'\uparrow T,x'\to x} \vr(t',x')= \hat g(x)$ for all $x\in \R$. If moreover there exists $\alpha\in (0,1)$ such that  $\hat g\in C^{4+\alpha}_{b}$,   $| \dxx \hat g| \leq \varepsilon^{-1}$ and $(T,\cdot, \dxx \hat g)\in \Dc\epu$ for some $\varepsilon>0$, then, for each $(t,x)\in [0,T)\x \R$, we can find $\phi \in \Ac$ such that $V^{t,x,v,\phi}_{T}=\hat g(X^{t,x,\phi}_{T})$
  with $v=\vr(t,x)$. 
 \end{theo}

 In \cite{BoLoZo2}, the authors also provide a viscosity solution characterization of  $\vr$, but in their case 
 \begin{itemize}
 \item[\rm (i)] admissible strategies should satisfy $\gamma\le \tilde \gamma(\cdot, X^{t,x,\phi})$ for some given function  $\tilde \gamma<\bar \gamma$ (uniformly on $[0,T]\x \R$),
 \item[\rm (ii)] $\bar F(\cdot, \tilde \gamma)<\infty$,
 \item[\rm (iii)] $\bar F(t,x,\cdot)$ is convex on $(-\infty,\tilde \gamma(t,x)]$ for all $(t,x)\in [0,T]\x \R$.
 \end{itemize}
 None of these assumptions are imposed here, and we also consider the case $\bar \gamma\equiv +\infty$.
 \vs2
 
 Still, the supersolution property can essentially be proved by mimicking the arguments of \cite[Section 5]{CheriSonTouz}, up to considering a weak formulation of our stochastic target problem. This will only provide a supersolution of \eqref{eq: pde vr} that will serve as a lower bound, see Proposition \ref{prop: visco sol vr eps} for a precise statement.   In \cite{BoLoZo2}, the subsolution property could not be proved directly as in  \cite{CheriSonTouz}. The reason is that the feedback effect of the controled state dynamics $(X,Y,V)$ prevented them to establish the required geometric dynamic programming principle. Instead, they used a smoothing argument in the spirit of \cite{krylov2000rate}. This however requires $\bar F$ to be convex, which, again, is not the case in our generalized setting.  We will instead rely on the theory of parabolic equations, which, up to regularization arguments, will allow us to construct smooth subsolutions of \eqref{eq: pde vr} from which superhedging strategies can be deduced, see Corollary \ref{cor: vr=u for Phi=g}. As in \cite{BoLoZo2}, combining these two results will prove Theorem \ref{theo: vr=u for Phi=hat g 0}. 
 \vs2

 We conclude this section with a remark on our definition of the face-lift of $g$. 
  \begin{remarque}\label{rem: face-lift} In \cite{BoLoZo2}, the face-lift  is defined as the smallest function above $g$ that is a viscosity supersolution of the equation $\bar \gamma - \partial^{2}_{x}\vp=  0$. It is obtained by considering any twice continuously differentiable function $\bar \Gamma$ such that $\partial^{2}_{x} \bar \Gamma=\bar \gamma$, and then setting
\begin{equation*}\label{eq: def hat g}
\bar{g}:=(g-\bar \Gamma)^{\text{conc}} +  \bar \Gamma,
\end{equation*}
in which the superscript \text{conc} means concave envelope, cf.~\cite[Lemma 3.1]{SonTouz}. This {actually} corresponds to our definition. The fact that $\hat g\ge \bar g$ is trivialy deduced from the supersolution property in the definition of $\bar g$. Let us prove the converse inequality. Fix $\eps\in (0,\eps_{\circ}]$, and define $\bar g_{\eps}$ as $\bar g$ but with $\bar \gamma-\eps$ in place of $\bar \gamma$. Fix  $\psi\in C^{\infty}_{b}$ with compact support, such that $\int \psi(y)dy=1$ and $\psi\ge 0$, and define $\bar g^{\eps}_{n}(x):=\int \bar g_{\eps}(y) \frac1n\psi(n({y-x}))dy$ for $n\ge 1$. Since $\bar g_{\eps}$ is the sum of a concave function and a $C^{2}$ function, one can consider the measure $m_{\eps}$ associated to its second derivative and it satisfies 
$m_{\eps}(dy)\le (\bar \gamma(y)-\eps) dy$. Then, $\partial^{2}_{x}\bar g^{\eps}_{n}(x)=\int \bar g_{\eps}(y) n \partial^{2}_{x}\psi(n({y-x}))dy$ $=$ $\int \frac1n \psi(n({y-x}))dm_{\eps}(y)$ $\le$ $\int \frac1n \psi(n({y-x}))(\bar \gamma(y)-\eps)dy$. Now, note  that $\bar g$ is continuous and therefore uniformly continuous on compact sets. Then, up to using the approximation from above argument of \cite[Lemma 3.2]{BoLoZo2}, we can assume that it is uniformly continuous. 
Since $\bar \gamma$ is also uniformly continuous, see \eqref{eq: hyp bar gamma unif cont}, one can find $\kappa,\eps>0$ such that $\bar g_{n}^{\eps,\kappa}: x\in \R\mapsto \bar g^{\eps}_{n}(x)+\kappa$ is $C^{2}$,   $\partial^{2}_{x}\bar g_{n}^{\eps,\kappa}\le \bar \gamma$ and $\bar g_{n}^{\eps,\kappa}\ge g$. By definition, it follows that $\bar g_{n}^{\eps,\kappa}\ge \hat g$. Clearly, $(\bar g_{n}^{\eps,\kappa})_{\eps,\kappa>0,n\ge 1}$ converges pointwise to $\bar g$ as $n\to \infty$ and $(\eps,\kappa)\to 0$ in a suitable way. This shows that $\bar g\ge \hat g$. 
  \end{remarque}

 \subsection{Supersolution property of a lower bound and partial comparison}
  
 In this section, we produce a supersolution of a version of \eqref{eq: pde vr} that is associated to $ \vr^{\eps}$, recall \eqref{eq: def v veps vepsk}, and that is a lower bound for  $\vr^{\eps}$. We also prove a partial comparison result on  this version that will be of important use later on.  Recall the definition of $\hat g^{\eps}$ in \eqref{eq: def hat g eps}.
 
 \begin{prop} \label{prop: visco sol vr eps}For each $\eps\in (0,\eps_{\circ}]$ small enough, there exists a continuous  function $\vrb^{\eps}\le \vr^{\eps}$ that has linear growth, is bounded from below,  is a viscosity super-solution of  
 \begin{align}\label{eq: pde vr eps}
 \min\{-\partial_{t}\vp- \bar F(\cdot,\dxx \vp)\;,\; \eps^{-1} - F(\cdot,\dxx \vp)\}=0\mbox{ on } [0,T)\x \R \tag{{\rm Eq}$_{\eps}$}
  \end{align}
 and satisfies $\liminf_{t'\uparrow T,x'\to x} \vrb^{\eps}(t',x')\ge \hat g^{\eps}(x)$ for all $x\in \R$. 
 \end{prop}
 \proof This follows from   exactly the same arguments as in \cite[Section 3.1]{BoLoZo2}.   We only explain the differences. As in  \cite[Section 3.2]{BoLoZo2}, we first introduce a sequence of weak formulations.  On  $(C(\R_{+}))^{5}$, let us   denote by $(\tilde \zeta:=(\tilde \gamma,\tilde b,\tilde \alpha,\tilde \beta), \tilde W)$ the coordinate process and let $\tilde {\mathbb F}^{\circ}=(\tilde \Fc^{\circ}_{s})_{  s \le T}$ be its raw filtration.  We say that a probability measure $\tilde \bP$ belongs to $\tilde \Ac_{k}$ if $\tilde W$ is a $\tilde \bP$-Brownian motion and  if for all $0\le \delta\le 1$ and $r\ge 0$ it holds $\tilde \bP$-a.s.~that
\begin{align} \label{eq: cond weak 1}
\tilde \gamma=\tilde \gamma_{0}+\int_{0}^{\cdot}\tilde \beta_{s} ds +\int_{0}^{\cdot} \tilde \alpha_{s} d\tilde W_{s}\;\mbox{ for some } \tilde \gamma_{0}\in \R,
\end{align}
\begin{align} \label{eq: cond weak 2}
 \sup_{\R_{+}} |\tilde \zeta|\le k\;,\;
\end{align}
and
\begin{align} \label{eq: cond weak 3}
\E^{{\tilde \bP}}\left[\sup\left\{|\tilde \zeta_{s'}-\tilde \zeta_{s}|,\; r\le s \le s'\le s+\delta\right\}|{\tilde \Fc_{r}^{\circ}} \right]\le k\delta.
\end{align}
For $\tilde \phi:=(y,\tilde \gamma,\tilde b)$, $y\in \R$, we define $(\tilde X^{x,\tilde \phi},\tilde Y^{\tilde \phi}, \tilde V^{x,v,\tilde \phi})$ as in \eqref{eq: def X}-\eqref{eq: def Y}-\eqref{eq: def V} associated to the control $(\tilde \gamma,\tilde b)$ with time-$t$ initial condition $(x,y,v)$, and  with $\tilde W$ in place of $W$. {For $t\le T$ and $k\ge 1$, we} say that $\tilde \bP\in \tilde \Gc_{k,\eps}(t,x,v,y)$ if
\begin{equation} \label{eq: def sur rep sous contrainte faible}
\left[\tilde  V^{x,v,\tilde \phi}_{{ T}}\ge g(\tilde X^{x,\tilde \phi}_{{ T}}),\; F(\cdot,\tilde X^{x,\tilde \phi},\tilde \gamma)\le \eps^{-1} \;\mbox{ and } \;   \tilde \gamma\in [-k,k] \mbox{ on } \R_{+}\right]\;\;\tilde \bP-{\rm a.s.}
 \end{equation}
We finally define
\begin{equation*}\label{eq: def v weak k}
   \vrt^{\eps}_{k}(t,x):=\inf\{ v=c+ yx~:~(c,y)\in \R\times [-k,k] \mbox{ s.t. }   \tilde \Ac_{k}\cap   \tilde \Gc_{k,\eps}(t,x,v,y)\ne \emptyset\}.
 \end{equation*}
 {\sl Step 1.} We first provide bounds for $\vrt^{\eps}_{k}$.  Note that $\vrt^{\eps}_{k}\le \vr^{\eps}_{k}$, so that \eqref{eq: bound vr et vr eps k} implies that  $[\vrt^{\eps}_{k}]^{+}$ has linear growth, uniformly in  $k\ge k_{\circ}$. Moreover, note that the fact that $\sigma$ is Lipschitz with linear growth in its second variable, uniformly on    $\Dc_{\eps,k}\x \R$ (see \eqref{eq: hyp sigma}), implies that 
  $\tilde X^{t,x,\tilde\phi}$ is a square integrable martingale under $\tilde \bP$ for any $\tilde \phi:=(y,\tilde \gamma,\tilde b)$, and that the same holds for $\int_{t}^{\cdot} \tilde Y^{t,\tilde\phi}_{s} d\tilde X^{t,x,\tilde\phi}_{s}$. Then, the inequality 
  $$
  v+\int_{t}^{T} F(s,\tilde X^{t,x,\tilde\phi}_{s},\tilde \gamma_{s})ds + \int_{t}^{T} \tilde Y^{t,\tilde\phi}_{s} d\tilde X^{t,x,\tilde\phi}_{s}\ge g(\tilde X^{t,x,\tilde\phi}_{T})
  $$
  combined with \eqref{eq: bound F sur Dc eps} and \eqref{eq: hyp g} implies that $v\ge -\sup|g^{-}|- T \sup_{\Dc_{\eps}} F>-\infty$. This shows that $\vrt^{\eps}_{k}$ is bounded from below, uniformly in $k\ge k_{0}$. 
  
  {\sl Step 2.} We claim that 
  \begin{equation*}
   \vrb^{\eps}(t,x):=\liminf_{\tiny \begin{array}{c}(k,t',x')\to (\infty,t,x)\\(t',x')\in [0,T)\x \R\end{array}}\vrt^{\eps}_{k}(t',x'),\;\;\;(t,x)\in [0,T]\x \R,
 \end{equation*}
  is a  viscosity  supersolution of \eqref{eq: pde vr eps}. 
   To prove this, it  suffices to show that it holds for each $\vrt^{\eps}_{k}$, with $k\ge k_{\circ}$,  and then to apply standard stability results, see e.g.~\cite{Bar94}. By the same arguments as in \cite[Proposition 3.15]{BoLoZo2}, each $\vrt^{\eps}_{k}$ is lower-semicontinuous.   Given a  $C^{\infty}_{b}$ test function $\vp$ and  $(t_{0},x_{0})\in [0,T)\x \R$ such that
$$
\text{(strict)}\min_{[0,T)\times\mathbb{R}} (\vrt^{\eps}_{k}-\vp)=(\vrt^{\eps}_{k}-\vp)(t_0, x_0)=0,
$$
we first use \eqref{eq: bound inf sigma} and the arguments of \cite[Step 1-2, proof of Theorem 3.16]{BoLoZo2} to obtain that there exists $\tilde \gamma_{0}$ such that 
$$
\dxx \vp(t_{0},x_{0})\le \tilde \gamma_{0}\mbox{ and } F(t_{0},x_{0},\tilde \gamma_{0})\le \eps^{-1}. 
$$
Then, the same arguments as in \cite[Step 3.a., proof of Theorem 3.16]{BoLoZo2} combined with \eqref{eq: def bar F} and  \eqref{eq: bar F unif elliptic} 
lead to
\begin{align*}
0\le&  F(t_{0},x_{0},\tilde \gamma_{0})-\partial_{t} \vp(t_{0},x_{0})-\frac12 \sigma^{2}(t_{0},x_{0},\tilde \gamma_{0})^{2} \dxx \vp(t_{0},x_{0})\\
&-\frac12\left(\tilde \gamma_{0}-\dxx \vp(t_{0}, x_{0})\right)\sigma^{2}(t_{0},x_{0},\tilde \gamma_{0})\\
=&-\partial_{t} \vp(t_{0},x_{0})-\bar F(t_{0},x_{0},\tilde \gamma_{0})\\
\le& -\partial_{t} \vp(t_{0},x_{0})-\bar F(t_{0},x_{0},\dxx \vp(t_{0}, x_{0})). 
\end{align*} 
Finally, the $T$-boundary condition is obtained as in \cite[Step 3.b., proof of Theorem 3.16]{BoLoZo2}, recall our assumption \eqref{eq: hyp g}, as well as Remark \ref{rem: face-lift}.
 \endp\\

 We now provide a partial comparison result that will be used later on. Note that a full comparison result could be proved as in \cite[Theorem 3.11]{BoLoZo2} when $\bar F$ is convex, by mimicking their arguments. It is however  not the case in general. Given the strategy of our proof, it is not required in this paper.  {In the following, we interpret  {\rm (Eq$_{\eps}$)} by using the convention $0^{-1}=\infty$ in the case $\eps=0$.}
 
 \begin{prop}\label{prop: Comp}  Let $U$   be an upper semicontinuous viscosity subsolution  of \eqref{eq: pde vr eps} for  $\eps\in [0,\eps_{\circ}]$. Let $V$ be a lower semicontinuous viscosity supersolution   of  {\rm (Eq$_{\eps'}$)}   for some $\eps'\in(\eps,\eps_{\circ}]$.  Assume that $U$ and $V$ have linear growth and that  $U\le V$ on $\{T\}\x \R$, then   $U\leq V$ on $[0,T] \times \mathbb{R}$.
\end{prop}
 
  \proof  Set $\hat U(t,x):=e^{\rho t}U(t,x)$, $\hat V(t,x):=e^{\rho t}V(t,x)$ for some $\rho>0$. Then, $\hat U$ is a subsolution of 
\begin{align} \label{eq: replacedComp1}
 \min\left\{\rho\vp-\partial_{t}\vp- e^{\rho \cdot}  \bar F(\cdot, e^{-\rho \cdot} \dxx \vp), \eps^{-1}-F(\cdot,e^{-\rho \cdot} \dxx \vp)\right\}=0
\end{align}   
   and $\hat V$  is a supersolution of
\begin{align} \label{eq: replacedComp2}
 \min\left\{\rho\vp-\partial_{t}\vp-e^{\rho \cdot} \bar F(\cdot, e^{-\rho \cdot} \dxx \vp),(\eps')^{-1}-F(\cdot,e^{-\rho \cdot} \dxx \vp)\right\}=0
\end{align}
on $ [0,T) \times \mathbb{R}.$ 

Assume that  $\sup_{[0,T]\times\mathbb{R}}(\hat U-\hat V)>0$. 
Then, there exists $\eta>0$ such that, {for all $n>0$}  and all $\lambda>0$ small enough,
\begin{equation}\label{eq: Thetan}
 \sup_{(t,x,y)\in [0,T]\times\mathbb{R}^{2}}\Big[\hat U(t,x)-\hat V(t,y)-\frac{\lambda}{2}|x|^{2}- \frac{n}{2} |x-y|^{2} \Big]\ge \eta>0.
\end{equation}
 Denote by  $(t_{n},x_{n},y_{n})$ the point at which this supremum is achieved. Since $\hat V(T,\cdot)\ge \hat U(T,\cdot)$, we have $t_{n}<T$. Moreover, 
  standard arguments, see e.g.,~\cite[Proposition 3.7]{CrandallIshiiLions}, lead to 
\begin{eqnarray}\label{eq: dif xn yn}
\lim\limits_{n \to \infty} n|x_{n}-y_{n}|^{2}=0.
\end{eqnarray}

We now apply  Ishii's lemma to obtain  the existence of $(a_{n},M_{n}, N_{n})\in \R^{3}$ such that
\begin{eqnarray*}
&\left(a_{n}, n(x_{n}-y_{n}){+\lambda x_{n}}, M_{n}\right)\in \bar {\cal P}^{2,+}\hat U(t_{n},x_{n})\\
&\left(a_{n},-n(x_{n}-y_{n}) , N_{n} \right)\in \bar {\cal P}^{2,-}\hat V(t_{n},y_{n}),
\end{eqnarray*}
in which $\bar{\cal P}^{2,+}$ and  $\bar{\cal P}^{2,-}$ denote as usual the {\sl closed} parabolic super- and subjets, see \cite{CrandallIshiiLions}, and
\begin{eqnarray*}
\begin{pmatrix}
M_{n} & 0 \\
0 & -N_{n}
\end{pmatrix}
\leq
3n
\begin{pmatrix}
1 & -1 \\
-1 & 1
\end{pmatrix}
+
\begin{pmatrix}
3\lambda+\frac{{\lambda^{2}}}{n} & -\lambda \\
-\lambda & 0
\end{pmatrix}.
\end{eqnarray*}
In particular,
$M_{n}\le N_{n}+{2}\lambda+\lambda^{2}/n.$
Since $\hat V$ is a supersolution of \eqref{eq: replacedComp2} and $\eps< \eps'$, \eqref{eq: hyp F monotone unif cont} and \eqref{eq: dif xn yn} imply that $F(t_{n},x_{n},e^{-\rho t_{n}}M_{n})<\eps^{-1}$ for $\lambda>0$ small enough and $n$ large enough.  Hence, 
\begin{align*}  
 \rho\hat U(t_{n},x_{n})  -a_{n}-e^{\rho t_{n}} \bar F(t_{n},x_{n}, e^{-\rho t_{n}} M_{n})\le 0. 
\end{align*} 
On the other hand, the supersolution property of $\hat V$ combined with {\eqref{eq: bar F uniform cont sur Deps}} 
and 
\eqref{eq: bar F unif elliptic} 
implies that 
\begin{align*}  
0\le& \rho\hat V(t_{n},y_{n})  -a_{n}-e^{\rho t_{n}} \bar F(t_{n},y_{n}, e^{-\rho t_{n}} N_{n})\\ 
\le&  \rho\hat V(t_{n},y_{n})  -a_{n}-e^{\rho t_{n}} \bar F(t_{n},y_{n}, e^{-\rho t_{n}} M_{n})+e^{\rho t_{n}} \delta(e^{-\rho t_{n}} ({2}\lambda +\lambda^{2}/n))
\end{align*} 
in which $\delta(z)\to 0$ as $z\to 0$. Hence, 
\begin{align*}  
\rho(\hat U(t_{n},x_{n})-\hat V(t_{n},y_{n})) \le& \; e^{\rho t_{n}} \left(\bar F(t_{n},x_{n}, e^{-\rho t_{n}}M_{n})- \bar F(t_{n},y_{n}, e^{-\rho t_{n}} M_{n})\right) \\
&+e^{\rho t_{n}}\delta(e^{-\rho t_{n}} (\lambda +\lambda^{2}/n)).
\end{align*} 
Recalling \eqref{eq: dif xn yn} and {\eqref{eq: bar F uniform cont sur Deps}}, 
we obtain a contradiction to \eqref{eq: Thetan} by sending $n\to \infty$ and then $\lambda \to 0$.
 \endp

 
 \subsection{Regularity of solutions to \eqref{eq: pde vr eps}}\label{sec: regul}
 
 To complete the characterization of Proposition \ref{prop: visco sol vr eps}, we now study the regularity of solutions to \eqref{eq: pde vr eps}. We shall indeed show that \eqref{eq: pde vr eps} admits a smooth solution  $u^{\eps}$  such that $(\cdot,\dxx u)\in  \Dc_{\eps}$ on $[0,T]\x \R$, for   $\eps>0$ small enough  and for a certain class of terminal conditions. A simple verification argument will then show that $u^{\eps}$ dominates the super-hedging price $\vr$ if the terminal data $\Phi^{\eps}$ associated to $u^{\eps}$ dominates $\hat g$. A lower   bound $u_{\eps}$ for $\vr$ can also be constructed   by considering a terminal condition $\Phi_{\eps}\le \hat g$ and using our comparison result of Proposition \ref{prop: Comp} combined with  Proposition \ref{prop: visco sol vr eps}. Then, letting $\Phi_{\eps},\Phi^{\eps}\to \hat g$ in a suitable way  will be enough to show that $\vr$ is actually a  solution of {\rm (Eq$_{0}$)}, i.e.~to conclude the proof of Theorem   \ref{theo: vr=u for Phi=hat g 0}. 
 
The strategy we employ consists in  establishing {\rm a priori} estimates for the second derivative of the solution to  \eqref{eq: pde vr eps}. Once established, the equation becomes uniformly parabolic, and higher regularity follows  by standard parabolic regularity (see \cite{Lieberman}). Then, the continuity method (see \cite{GT}) allows us to actually construct the solution to \eqref{eq: pde vr eps}.

Let us start with uniform estimates for solutions to  \eqref{eq: pde vr eps} such that $(\cdot, \cdot,\dxx u) \in \Dc_{\eps'}$ for some $\eps'>0$, in the case where the terminal condition $\Phi$ is smooth and satisfies a similar constraint. 

  \begin{prop}\label{prop: u in DC bar gamma eps}  Let $u$ and $\Phi$   be two continuous functions such that  
\begin{enumerate}
\item[\rm (i)]  $\Phi \in C^{2}(\R)$ with $|\dxx \Phi| \leq K_\Phi$ for some $K_\Phi>0$,
\item[\rm (ii)] $(T,\cdot,\dxx \Phi) \in \Dc_{\eps_\Phi}$ for some $\eps_\Phi>0$,    
\item[\rm (iii)] $u\in C^{1,4}([0,T)\x \R) {\cap C^{0,2}}([0,T]\x \R)$ with $|\dxx u| \leq K'$ for some $K'>0$,
\item[\rm (iv)] $(\cdot, \cdot,\dxx u) \in \Dc_{\eps'}$ for some $\eps'>0$.
\end{enumerate}
Assume that $u$ solves  
 \begin{align}
 \partial_{t} u + \bar F(\cdot, \dxx u)&=0 \;\;\mbox{ on } [0,T)\x \R,\label{eq: pde u inside}\tag{{\rm Eq}$_{0}$}\\
 u(T,\cdot)&=\Phi \;\;\mbox{ on } \R.\label{eq: pde u T}
 \end{align}
 Then, 
 \begin{itemize}
\item[a.]   $(\cdot,\dxx u)\in  \Dc_{\eps}$ on $[0,T]\x \R$, for some $\eps>0$ that depends only on $\eps_\Phi$ and $L_{\circ}$, 
\item[b.] $|\dxx u| \leq K$ on $[0,T]\x \R$ where $K$ depends only on $K_\Phi$.
\item[c.] If  $\Phi$ is globally Lipschitz, then $u$ is also globally Lipschitz with Lipschitz constant controlled by the one of $\Phi$. 
\item[d.] $u$ is  the unique $C^{1,2}([0,T)\x \R)\cap C^{0}([0,T]\x \R)$ solution of \eqref{eq: pde u inside}-\eqref{eq: pde u T} such that  $(\cdot, \cdot,\dxx u) \in \Dc_{\eps''}$ for some $\eps''>0$.
\end{itemize}

 \end{prop} 
 
 \proof  a. Let $V:=\bar F(\cdot, \dxx u)$. Then, on $[0,T)\x \R$,
 \begin{align*}
 \partial_{t} V=\partial_{t}\bar F(\cdot, \dxx u)+\partial_{z}\bar F(\cdot, \dxx u)\partial_{t}\dxx u
 \end{align*}
 in which, by \eqref{eq: pde u inside}, $\partial_{t}\dxx u+\dxx V=0$. Hence, 
 \begin{align}\label{eq: pde V 1}
 \partial_{t} V+\partial_{z}\bar F(\cdot, \dxx u) \dxx V=\partial_{t}\bar F(\cdot, \dxx u)=\frac{\partial_{t}\bar F(\cdot, \dxx u)}{\bar F(\cdot, \dxx u)} V,
 \end{align}
recall \eqref{eq: bar Ft/ bar F bounded by Lo}. 
For $(t,x)\in [0,T]\x \R$, let  $\bar X^{t,x}$ be the solution of 
$$
\bar X=x+\int_{t}^{\cdot} (2\partial_{z}\bar F(\cdot, \dxx u)(s,\bar X_{s}))^{\frac12} dW_{s}. 
$$  
By (iv), \eqref{eq : regul bar F cas general}  and \eqref{eq: bar F unif elliptic}, it is well-defined. Combining It\^{o}'s Lemma and a standard localizing argument using  \eqref{eq : regul bar F cas general} 
and \eqref{eq: bar Ft/ bar F bounded by Lo}, we obtain 
\begin{align}\label{eq: Ito sur V up to T}
V(t,x)=\E[V(T,\bar X^{t,x}_{T})e^{-\int_{t}^{T} (\partial_{t}\bar F(\cdot, \dxx u)/\bar F(\cdot, \dxx u))(s,\bar X^{t,x}_{s}) ds}].
\end{align}
By definition of $V$ and the fact that $\dxx u(T,\cdot)=\dxx \Phi$ by (iii), this shows that $(\cdot,\dxx u)\in  \Dc_{\eps}$ on $[0,T]\x \R$, for some $\eps>0$ that depends only on $L_\circ$ and $\eps_\Phi$.

b. To obtain the bound on $\dxx u$, we first differentiate twice \eqref{eq: pde u inside} with respect to $x$, recall \eqref{eq : regul bar F cas general} and (iii).  Letting $Z(t,x)=\dxx u(t,x)$,  this yields
\ben
\dt  Z +  2\dx \partial_z \bar{F} \partial_x Z + \partial_z\bar{F}\dxx Z + \partial^2_z \bar{F}(\dx Z)^2=-\dxx \bar{F}.
\enn
We now consider $$(t,x)\mapsto\underline{Z}(t,x) := \min\{0, \inf Z(T,\cdot)\}e^{M(T-t)},$$
in which $M$ is given in \eqref{eq: bar Ft/ bar F bounded by Lo}. 
Then, 
\ben
\dt  \underline{Z} +  2\dx \partial_z \bar{F} \partial_x \underline{Z} + \partial_z\bar{F}\dxx \underline{Z} + \partial^2_z \bar{F}(\partial_{x}\underline{Z})^2 &=& -M\underline{Z}
\geq -\dxx \bar{F}(t,x,\underline{Z}).
\enn 
Under the current assumptions, $Z$ is uniformly bounded on $[0,T]\times\R$. Moreover, from assumption (\ref{eq : regul bar F cas general}), 
$\dxx \bar{F}$ is uniformly continuous on $\Dc_{\eps,\eps^{-1}}$, for all $\eps>0$ small enough, 
hence, by  \eqref{eq: bar F(,0)=0} and \cite[Proof of comparison, Theorem 5.1]{CrandallIshiiLions},  the  comparison principle holds between $Z$ and $\underline{Z}$, and yields that $\underline{Z} \leq Z$ globally on $[0,T]\times \R$. The upper bound is obtained in the exact same way.

c.  The assertion about the Lipschitz regularity   also follows from  the linearised equation satisfied by $\kappa=\dx u $:
\ben
\dt \kappa + \partial_{z}\bar F(\cdot,\dxx u) \dxx \kappa + \dx\bar F(\cdot,\dx\kappa)=0,\;\kappa(T,\cdot)=\dx \Phi.
\enn
Under the assumptions   \eqref{eq: bar F unif elliptic}, \eqref{eq : regul bar F cas general}, and \eqref{eq: bar F(,0)=0}, 
this implies that 
$$
\kappa(t,x)=\E[\dx \Phi(\tilde X^{t,x}_{T})]
$$
where 
$$
\tilde X^{t,x}=x+\int_{t}^{\cdot} \left(2 \partial_{z}\bar F(\cdot,\dxx u) \right)^{\frac12} {(s,\tilde X^{t,x}_{s})}dW_{s}+\int_t^\cdot \frac{\dx \bar F({\cdot,\dxx u})}{\dxx u}{(s,\tilde X^{t,x}_{s})}ds,
$$
and the result follows.
  (Note that, since $\bar F(\cdot,0)=0$ and  $\bar F  \in C^{1,3,3}_b(\Dc_{\eps,\eps^{-1}})$,  the map 
{$z\mapsto\frac{\dx \bar F(\cdot,z)}{z}$} is bounded and Lipschitz - after extending it to $\partial_{z}\dx \bar F(\cdot,0)$ at $0$.)

d. Consider another solution $u'$. Then, b.~implies that $u$ and $u' $ have at most a quadratic growth. Moreover, a.~allows one to consider a uniformly parabolic equation. Then, the fact that $u=u'$ follows from standard arguments.  
\endp
\vs2
 
We are now in position to construct a smooth solution to \eqref{eq: pde u inside}. 
\begin{theo}\label{theo: u C4}
Let $\Phi$ be a continuous map such that $| \dxx \Phi| \leq \varepsilon^{-1}$ and $(T,\cdot, \dxx \Phi)\in \Dc\epu$ for some $\varepsilon>0$. Then, there exists a solution $u$ of \eqref{eq: pde u inside}-\eqref{eq: pde u T} that belongs to $C([0,T]\times \R) \cap C^{1,4}_{loc}([0,T)\times \R)$, such that  $|\dxx u|\le (\eps_{\Phi,L_{\circ}})^{-1}$   and $(\cdot,\dxx u)\in \Dc_{\eps_{\Phi,L_{\circ}}}$  on $[0,T]\x \R$, for some $\eps_{\Phi,L_{\circ}}>0$ that only depends on $\Phi$ and $L_{\circ}$.  If  $\Phi$ is globally Lipschitz, then $u$ is also globally Lipschitz with Lipschitz constant controlled by the one of $\Phi$. If moreover there exists $\alpha\in (0,1)$ such that  $\Phi \in C^{4+\alpha}_b$ then $u \in C^{1,4}_{b}$.
\end{theo}

\proof  
This follows by using the continuity method (cf. \cite[Chap.~17.2]{GT}). We first mollify $\Phi$  into a function $\Phi_n$ so that $\partial^5_x \Phi_{n}$ is bounded, and at the same time $\bar F$ so that $\bar F(\cdot,\cdot,z)\in C^\infty([0,T]\times\R)$ locally uniformly with respect to $z$. This is possible, since  $\bar\gamma$ and $\bar F$ are uniformly continuous (recall \eqref{eq: hyp bar gamma unif cont} and \eqref{eq: bar F uniform cont sur Deps}), 
by taking a compactly supported smoothing kernel $\psi \in C^\infty(\R)$ and considering 
\ben
\Phi_n &=& \frac1n\int_{\R} \Phi(y) \psi(n({y-\cdot}))dy,\;\\
\bar F_n& =& \frac{1}{n^{2}}\int_{[0,T]\x \R} \bar F(s,y,\cdot){\psi(n(s-\cdot))\psi(n(y-\cdot))}dsdy.
\enn
For later use, note that $\bar F_{n}(T,\cdot,\partial^{2}_{x}\Phi_n)\le 2\eps^{-1}$, for $n$ large enough.
Set 
$$
G_{n}(\vp,\theta):=[\partial_{t} \vp + \bar F_n(\cdot, \partial^{2}_{x} \vp)]\1_{[0,T)}+\1_{\{T\}}(\vp-\theta \Phi_n)\;\mbox{ for $\vp \in C^{1,4}_b$,}
$$
and let  $E_{n}\subset [0,1]$ be the set of real number $\theta \in [0,1]$ for which   a  $C^{1,4}_b$ solution $u^{n}_\theta$  to $G_{n}(u^{n}_{\theta},\theta)=0$ exists such that it satisfies the condition (iii)-(iv) of Proposition \ref{prop: u in DC bar gamma eps}.
By \eqref{eq: bar F(,0)=0}, $u_{0}\equiv 0$ solves $G_{n}(u_{0},0)=0$ so that  $0\in E_{n}$. Hence, $E_{n}$ is non empty. 
Moreover, for every $\theta \in E_{n}$, the linearised operator associated to $G_{n}$ is
$$
 (\tilde u,\tilde \theta)\in C^{1,2}\x E_{n}\mapsto L_{n}(\tilde u,\tilde \theta):= [\dt \tilde u  +\partial_z \bar F_{n}(t,x, \dxx u) \dxx \tilde u]\1_{[0,T)}+\1_{\{T\}}(\tilde u-\tilde \theta \Phi_n).
 $$
It is uniformly parabolic (recall \eqref{eq: bar F unif elliptic}) with coefficients in $C^\infty$. For $\tilde \theta$ fixed, the equation $L_{n}(\tilde u,\tilde \theta)=0$   is therefore a linear, uniformly parabolic  equation, with smooth coefficients. The terminal data is smooth, has linear growth and bounded derivatives of order 1 up to 5.
 Standard parabolic regularity theory (see \cite{Friedman}) yields that the linearised  equation with respect to $u$ is solvable in $C^{1,4}_b$.
By the implicit function theorem, see e.g.~\cite[Theorem 17.6]{GT}, $E_{n}$ is  open  {in $[0,1]$}. By the a priori estimates of Proposition \ref{prop: u in DC bar gamma eps}, $E_{n}$ is also closed. Therefore, $E_{n}=[0,1]$ and $u^{n}_{1}$ is well defined. Since $(\bar F_{n})_{n\ge 1}$ is uniformly parabolic, uniformly in $n$, and $(\Phi_{n})_{n\ge 1}$ is bounded in $C^{4+\alpha}_{b}$ uniformly in $n$, then   $(u^{n}_{1})_{n\ge 1}$ is $C^{1,4}_{b}$ uniformly in $n$.   It remains to send $n\to \infty$ and to appeal again to the  a priori estimates of Proposition \ref{prop: u in DC bar gamma eps} to deduce the required result. \endp

 \subsection{Full chacterization of the super-hedging price and perfect hedging in the smooth case}
  
We are now about to conclude the proof of Theorem \ref{theo: vr=u for Phi=hat g 0}. Let 
$\hat u$ be the function  constructed in  Theorem \ref{theo: u C4}  for $\Phi=\hat g$, assuming that $\hat g$ satisfies the required constraints.  We first establish that $\hat u$ permits to apply a perfect hedging strategy of the face-lifted payoff whenever it is smooth enough, and that it coincides with the super-hedging price.  
  
 \begin{cor}\label{cor: vr=u for Phi=g} Assume that there exists $\alpha\in (0,1)$ such that  $\hat g\in C^{4+\alpha}_{b}$, that $| \dxx \hat g| \leq \varepsilon^{-1}$ and $(T,\cdot, \dxx \hat g)\in \Dc\epu$ for some $\varepsilon>0$.  Let $\hat u$ be the function  constructed in  Theorem \ref{theo: u C4}  for $\Phi=\hat g$. Then, $\vr=\hat u$  and, for each $(t,x)\in [0,T]\x \R$, we can find $\phi \in \Ac$ such that $V^{t,x,v,\phi}_{T}=\hat g(X^{t,x,\phi}_{T})$.
 \end{cor}
 
 \proof  It follows from Theorem \ref{theo: u C4}, It\^{o}'s lemma and \eqref{eq: def bar F} that  $\hat u$  induces an exact replication strategy:
 \begin{align*}
 \hat g(X^{t,x,\phi}_{T})
 =& \hat u(t,x)+\int_{t}^{T} \left[\partial_{t} \hat u+ \frac12 \sigma(\cdot,\dxx \hat u)^{2} \dxx \hat u\right](s,X^{t,x,\phi}_{s})ds\\
 & +\int_{t}^{T} \dx \hat u(s,X^{t,x,\phi}_{s})dX^{t,x,\phi}_{s}\\
  =& \hat u(t,x)+\int_{t}^{T} F(s,X^{t,x,\phi}_{s},\gamma_{s})ds +\int_{t}^{T} Y^{t,x,\phi}_{s}dX^{t,x,\phi}_{s}
 \end{align*}
 in which 
 $\phi=(y,b,\gamma)$ with 
 $$
 y= \dx \hat u(t,x),\;b=  ([\partial_{t} + \frac12 \sigma(\cdot,\gamma)^{2} \dxx  ]\dx \hat u) (\cdot,X^{t,x,\phi}_{\cdot}) ,\;\gamma=\dxx \hat u(\cdot,X^{t,x,\phi}_{\cdot}).
 $$
 Hence, $\hat u\ge \vr$. Moreover, $\hat u$ is a viscosity subsolution of {\rm (Eq$_{\eps'}$)} for all $\eps'\ge 0$ small enough. Since $\hat g$ is globally  Lipschitz, $\hat u$ is also globally Lipschitz (Theorem \ref{theo: u C4}), and therefore has linear growth. By Proposition \ref{prop: visco sol vr eps}, $\vr^{\eps}\ge \vrb^{\eps}$ that is a super-solution of   \eqref{eq: pde vr eps} and satisfies $\liminf_{t'\uparrow T,x'\to x} \vrb^{\eps}(t',x')\ge \hat g^{\eps}(x)\ge \hat g(x)=\hat u(T,x)$ for all $x\in \R$.
 Then, Proposition \ref{prop: Comp} implies that $\vr^{\eps}\ge \hat u$. Taking the inf over $\eps>0$ leads to $\vr \ge \hat u$. 
 \endp
 \vs2
 
 We can now conclude the proof of Theorem \ref{theo: vr=u for Phi=hat g 0}.
 
 \noindent{\bf Proof of Theorem \ref{theo: vr=u for Phi=hat g 0}.}
 For  $\varepsilon>0$, let $\Phi\epu,\Phi\ep\in C^{2}$ be  such that, for $\Psi\in \{\Phi\epu,\Phi\ep\}$,   
\ben
\Psi \in C^{5}_{b}(\R),\;
| \partial^2_{x}\Psi |\leq \varepsilon^{-1} ,\;
(T,\cdot, \partial^2_{x}\Psi)\in \Dc_\eps,
\enn
and
\ben
\Phi\epu \leq \hat g \leq \Phi\ep,\;
\Phi\ep-\Phi\epu \leq \delta(\varepsilon),
\enn
in which  $\lim_{\varepsilon\to 0}\delta(\varepsilon)=0$. Such functions can be constructed as in Remark \ref{rem: face-lift}, and we can further assume that $\Phi\ep$ (resp. $\Phi\epu$) is non-increasing (resp. non-decreasing) with respect to $\eps$.
Let $u\ep$ and $u\epu$ be the (smooth) solutions to {\rm (Eq$_{0}$)} associated to $\Phi\ep$ and $\Phi\epu$ respectively, as in Theorem \ref{theo: u C4}.
By applying Corollay \ref{cor: vr=u for Phi=g} to $\Phi\ep$ in place of $\hat g$, we deduce  that $u\ep$ is the super-hedging price of  $\Phi\ep\ge \hat g$ so that $u\ep \ge \vr$. Similarly $u\epu\le \vr$, and therefore  $u\epu \leq \vr \leq u\ep$.

By the comparison principle, we also have $$0\leq u\ep-  u\epu\le \sup \{\Phi\ep-\Phi\epu\}\le \delta(\eps).$$

It follows that  $\vr$ is the uniform limit of a sequence of continuous functions, and is therefore continuous. Each of the functions $u_{\eps}$ solves \eqref{eq: pde vr}, recall \eqref{eq:defbargamma}.  Standard stability results, see e.g.~\cite{Bar94}, imply that $\vr$ is  a viscosity solution to \eqref{eq: pde vr}-\eqref{eq: pde vr T}. 

The other assertions in Theorem \ref{theo: vr=u for Phi=hat g 0} are immediate consequences of Corollary \ref{cor: vr=u for Phi=g}.
\endproof

\section{Asymptotic analysis}\label{sec: asympto}

We now consider the case where the impact of the $\gamma$ process in the dynamics of $(X,V)$ is small. Our aim is to obtain an asymptotic expansion around an impact free model. More precisely, we consider the dynamics  
 \begin{align*}
X^{\epsilon,t,x,\phi}&=x+ {  \int_{t}^{\cdot}  \mu(s,X^{\epsilon,t,x,\phi}_{s}, \epsilon\gamma_{s}, \epsilon b_{s})ds }+  \int_{t}^{\cdot}  \sigma(s,X^{\epsilon,t,x,\phi}_{s},\epsilon \gamma_{s}) dW_{s}\\
V^{\epsilon,t,x,v,\phi}&=v+\int_{t}^{\cdot} \epsilon^{-1}F(s,X^{\epsilon,t,x,\phi}_{s},\epsilon \gamma_{s})ds + \int_{t}^{\cdot} Y^{\epsilon,t,x,\phi}_{s} dX^{\epsilon,t,x,\phi}_{s},\;\epsilon>0, 
\end{align*}
 and denote by $\vr^{\epsilon}$ the corresponding super-hedging price.

 We place ourself in the context of  Corollary \ref{cor: vr=u for Phi=g}  {for the coefficients $\mu(\cdot,{\epsilon}\cdot,{\epsilon}\cdot),$ $\sigma(\cdot,\epsilon\cdot)$ and $\epsilon^{-1}F(\cdot,\epsilon\cdot )$}. In particular, we assume that $\hat g \in C^{2}$ is such that $\epsilon^{-1}\bar F(T,\cdot,\epsilon \dxx \hat g)$ is bounded on $\R$, for $\epsilon>0$ small enough.
 \vs2
 
 In the following, we use the notation 
 $$
 (\bar F_{0},\partial^{n}_{z} \bar F_{0}):=(\bar F(\cdot,0), \partial^{n}_{z}\bar F(\cdot,0)), \;\mbox{ for $n= 1,2$}.
 $$

\begin{remarque}\label{rem : scaling sur exemple} Note that the model of  \cite{BoLoZo2} corresponds to  
\begin{align*}
\sigma(t,x,\epsilon z)=\frac{\sigma_\circ(t,x)}{1-\epsilon f(x)z}\;,\;
\epsilon^{-1}F(t,x,\epsilon z)=\frac12 \left(\frac{\sigma_\circ(t,x) z}{1-\epsilon f(x)z}\right)^{2} \epsilon f(x).
\end{align*}
Our scaling therefore amounts to consider a small impact function $x\mapsto \epsilon f(x)$. In order to interpret the result of Proposition \ref{prop: asymptotique} below, also observe that  
\begin{align*}
(2\partial_{z} \bar F_{0}(t,x))^{\frac12} =\sigma_\circ(t,x)  &\;\mbox{ and }\;
\partial^2_z  \bar F_{0}(t,x)=\sigma_\circ^{2}(t,x)f(x).
\end{align*}
\end{remarque}
Our expansion is performed around the solution $\vr^{0}$ of 
\begin{align}\label{eq: pde vr 0}
 \partial_{t} \vr^{0} + \partial_{z} \bar F_{0} \partial_x ^{2}\vr^{0}= 0 \;\mbox{ on $[0,T)\x \R$} \;\mbox{ and } \vr^{0}(T,\cdot)=\hat g \mbox{ on } \R.
\end{align} 

\begin{remarque}\label{rem: vr 0 est smooth} Let the conditions of Corollary \ref{cor: vr=u for Phi=g} hold and assume that $\bar F\in C^{1,3,1}_{loc}( \Dc)$ with 
\begin{align}\label{eq: hyp pour vr 0}
|\partial_{x}\partial_{z} \bar F_{0}|+ |\partial^{2}_{x}\partial_{z} \bar F_{0}|  \mbox{ uniformly bounded.}
\end{align}
Then, $\vr^{0}$ is the unique solution in $C^{1,2}_{b}([0,T]\x \R)\cap C^{1,3}([0,T)\x \R])$ of \eqref{eq: pde vr 0}. This follows from \eqref{eq: bar F unif elliptic} and standard estimates.
\end{remarque}

 The following expansion requires some additional regularity on $\hat g$ that will in general not be satisfied in applications. However, one can reduce to it up to a slight approximation argument. 
 \begin{prop}\label{prop: asymptotique} Assume that the conditions of Corollary \ref{cor: vr=u for Phi=g} hold  with $\bar F^{\epsilon}:=\epsilon^{-1}\bar F(\cdot,\epsilon \cdot)$ in place of $\bar F$, uniformly in $\epsilon \in (0,\epsilon_{\circ}]$, for some $\epsilon_{\circ}>0$. Assume further that   $\bar F\in C^{1,2,3}_{loc}( \Dc)$, that \eqref{eq: hyp pour vr 0} and   
 \begin{align}\label{eq: borne Fzz et Fzzz sur Deps}
  \sup_{\Dc_{\epsilon}}\left(|{\partial^{2}_{z}} \bar F_{0}|+|{\partial^{3}_{z}} \bar F_{0}|+ |\partial_{x}\partial^{2}_{z} \bar F_{0}|+ |\partial^{2}_{x}\partial^{2}_{z} \bar F_{0}|\right)<\infty
  \end{align}
  hold. 
Then, there exists some $o(\eps)$, which does not depend on $x$, such that 
\begin{align*}
\vr^{\epsilon}(0,x)=&\vr^{0}(0,x)+\frac{\epsilon}2\E\left[  \int_{0}^{T}   [\partial^2_z  \bar F_0 (\partial_x ^{2}\vr^{0})^{2}](s,\tilde X^{0}_{s})ds \right]+o(\epsilon)\\
=&\vr^{0}(0,x)+{\frac{\epsilon}2}\;\E\left[\partial_x \hat g(T,\tilde X^{0}_{T})\tilde Y_{T}\right]+o(\epsilon)
\end{align*}
where, for $z\in \R$,  $\tilde X^{z}$ is the solution on $[0,T]$ of 
  \begin{align}\label{eq: tilde Xz}
  \tilde X^{z}=x+\int_{t}^{\cdot} (2 \partial_{z} \bar F (\cdot,z \partial_x ^{2}\vr^{0}(\cdot)))^{\frac12}(s,\tilde X^{z}_{s}) dW_{s},
  \end{align}
and $\tilde Y:=\partial_{z} \tilde X^{z}{|_{z=0}}$, solves
  $$
  \tilde Y= \frac1{\sqrt2} \int_{t}^{\cdot}  \frac{\partial_x  \partial_{z} \bar F_{0}(s,\tilde X^{0}_{s})\tilde Y_{s}+\partial^2_z  \bar F_{0}\partial_x ^{2}\vr^{0}(s,\tilde X^{0}_{s})  }{\sqrt{\partial_{z} \bar F_{0}(s,\tilde X^{0}_{s})}}dW_{s}.
  $$
    \end{prop}
  
 \proof By Corollary \ref{cor: vr=u for Phi=g}, each $\vr^{\epsilon}$ associated to $\epsilon \in (0,\epsilon_{\circ}]$ solves 
 $$
 \partial_{t} \vr^{\epsilon} + \epsilon^{-1}\bar F(\cdot, \epsilon \partial_x ^{2} \vr^{\epsilon})=0.
 $$
 Moreover, it follows from our assumptions and Corollary \ref{cor: vr=u for Phi=g}  that    $(\cdot,\vr^{\epsilon})\in \Dc_{\epsilon}$ for all $\epsilon \in (0,\epsilon_{\circ}]$. Then,  the fact that $\bar F(\cdot,0)=0$  implies that  
 $$
 \partial_{t} \vr^{\epsilon} + \partial_{z} \bar F_{0} \partial_x ^{2} \vr^{\epsilon}+\frac12 \epsilon \partial^2_z  \bar F_{0} (\partial_x ^{2}\vr^{\epsilon})^{2} =O(\epsilon^{2}),
 $$
 in which  the $O(\epsilon^{2})$ is uniform since $|{\partial^{3}_{z}} \bar F_{0}|$ is uniformly bounded on $\Dc_{\eps}$ by assumption. Let $\Delta v^{\epsilon}:=(\vr^{\epsilon}-\vr^{0})/\epsilon$. By the above, \eqref{eq: pde vr 0} and Remark \ref{rem: vr 0 est smooth}, it solves 
\begin{align*}
O(\epsilon)=&
 \partial_{t} \Delta v^{\epsilon} + \partial_{z} \bar F_{0} \partial_x ^{2} \Delta v^{\epsilon}+\frac12  \partial^2_z  \bar F_{0} (\partial_x ^{2}\vr^{0})^{2}\\
 &+\frac12 \epsilon^{2} \partial^2_z  \bar F_{0} (\partial_x ^{2}\Delta v^{\epsilon})^{2}+\epsilon \partial^2_z  \bar F_{0} \partial_x ^{2}\Delta v^{\epsilon} \partial_x ^{2}\vr^{0},
\end{align*} 
in which $O(\epsilon)$ is uniform on $[0,T)\x \R$. 
By Theorem \ref{theo: u C4}, Remark \ref{rem: vr 0 est smooth}, and the same arguments as in this remark,   $(\partial_x ^{2}\Delta v^{\epsilon},$ $\partial^2_z  \bar F_{0},$ $\partial_x ^{2}\vr^{0} )_{0<\epsilon\le \epsilon_{\circ}}$ is  {locally} bounded.  
Since $\Delta v^{\epsilon}(T,\cdot)=0$, it follows that 
$$
\Delta v^{\epsilon}(0,x)=\E\left[  \frac12\int_{0}^{T}   [\partial^2_z  \bar F_0 (\partial_x ^{2}\vr^{0})^{2}](s,\tilde X^{0}_{s})ds \right]+O(\epsilon).
$$
Hence,  
$\Delta v:=\lim_{\epsilon \to 0} \Delta v^{\epsilon}$ is given by 
\begin{equation}\label{eq: Delta vr}
\Delta v(0,x)=\E\left[  \frac12\int_{0}^{T}   [\partial^2_z  \bar F_0 (\partial_x ^{2}\vr^{0})^{2}](s,\tilde X^{0}_{s})ds \right].
\end{equation}
Moreover,  $\partial_x \vr^{0}$ satisfies 
\begin{align}\label{eq: pde dx vr0}
 \partial_{t} (\partial_x \vr^{0}) + \partial_x \partial_{z} \bar F_{0} \partial_x ^{2} \vr^{0}+\partial_{z} \bar F_{0} \partial_x ^{2} (\partial_x \vr^{0})=0,
\end{align}
recall Remark \ref{rem: vr 0 est smooth}. 

Applying It\^{o}'s lemma to $\partial_x \vr^{0}(t,\tilde X^{0}_{t})\tilde Y_{t}$, we obtain
\begin{align*}
&d(\partial_x \vr^{0}(t,\tilde X^{0}_{t})\tilde Y_{t})=\partial_t\partial_x \vr^{0}(t,\tilde X^{0}_{t})\tilde Y_{t}dt+\partial^2_x \vr^{0}(t,\tilde X^{0}_{t})\tilde Y_{t}d\tilde X^{0}_{t}+\partial_x \vr^{0}(t,\tilde X^{0}_{t})d\tilde Y_{t}\\
&+\partial^2_x \vr^{0}(t,\tilde X^{0}_{t})d\langle\tilde Y,\tilde X^{0}\rangle_{t}+\frac{1}{2}\partial^2_x (\partial_x \vr^{0}(t,\tilde X^{0}_{t}))\tilde Y_{t}d\langle\tilde X^{0}\rangle_{t}\\
&=\left(\partial_t\partial_x \vr^{0}(t,\tilde X^{0}_{t})+\partial^2_x \vr^{0}(t,\tilde X^{0}_{t})\partial_x  \partial_{z} \bar F_{0}(t,\tilde X^{0}_{t})+\partial^2_x (\partial_x \vr^{0}(t,\tilde X^{0}_{t}))\partial_{z} \bar F_0 (t,\tilde X^{0}_{t})\right)\tilde Y_{t}dt\\
&+\partial^2_z  \bar F_{0}(t,\tilde X^{0}_{t})(\partial_x ^{2}\vr^{0}(t,\tilde X^{0}_{t}))^2dt+\partial^2_x \vr^{0}(t,\tilde X^{0}_{t})\tilde Y_{t}d\tilde X^{0}_{t}+\partial_x \vr^{0}(t,\tilde X^{0}_{t})d\tilde Y_{t}\\
&=\partial^2_z  \bar F_{0}(t,\tilde X^{0}_{t})(\partial_x ^{2}\vr^{0}(t,\tilde X^{0}_{t}))^2dt+\partial^2_x \vr^{0}(t,\tilde X^{0}_{t})\tilde Y_{t}d\tilde X^{0}_{t}+\partial_x \vr^{0}(t,\tilde X^{0}_{t})d\tilde Y_{t}
\end{align*}
where we use (\ref{eq: pde dx vr0}) to get the last equality.

Therefore, taking expectation on both sides, we have 
$$\E\left[\partial_x \vr^{0}(T,\tilde X^{0}_{T})\tilde Y_{T}\right]=\E\left[ \int_{0}^{T}   [\partial^2_z  \bar F_0 (\partial_x ^{2}\vr^{0})^{2}](s,\tilde X^{0}_{s})ds \right],$$ 
which leads to
$$
\Delta v(0,x)=\frac12\E\left[ \partial_x \vr^{0}(T,\tilde X^{0}_{T})\tilde Y_{T}\right]=\frac12\E\left[ \partial_x \hat g(T,\tilde X^{0}_{T})\tilde Y_{T}\right]. 
$$ 
 \endp
 
  \begin{remarque}\label{rem: edp Delta v}  For later use, note that the above proof implies that $\Delta v$ defined in \eqref{eq: Delta vr} satisfies 
  \begin{align*}
 \partial_{t} \Delta v  + \partial_{z} \bar F_{0} \partial_x ^{2} \Delta v +\frac12  \partial^2_z  \bar F_{0} (\partial_x ^{2}\vr^{0})^{2}=0\;\mbox{ on $[0,T)\x \R$.}
 \end{align*} 
 \end{remarque} 
 
  \begin{remarque} A more tractable formulation can be obtained in the particular case where $(\partial_z  \bar F_{0},\partial^2_z  \bar F_{0})=(\lambda_{1},\lambda_{2})$ is  constant and $\dx \partial_z  \bar F_{0}=0$. This is the case in the model of  \cite{BoLoZo2}, see Example \ref{exemple: BoLoZo},  whenever $\sigma_\circ$ and $f$ are constant, see e.g.~Remark \ref{rem : scaling sur exemple}. Then, $\partial_x \vr^{0}(\cdot,\tilde X^{0})=\partial_x \vr^{0}(0,x)+\int_{0}^{\cdot} \sqrt{2\lambda_{1}}\partial_x ^{2}\vr^{0}(s,\tilde X^{0}_{s})dW_{s}$ by \eqref{eq: pde dx vr0}, so that
 \begin{align*}
 \frac{\epsilon}2\E\left[  \int_{0}^{T}   [\partial^2_z  \bar F_0 (\partial_x ^{2}\vr^{0})^{2}](s,\tilde X^{0}_{s})ds \right]
 &= \frac{\epsilon\lambda_{2}}{4\lambda_{1}}\E\left[  \int_{0}^{T}   [\sqrt{2\lambda_{1}}  \partial_x ^{2}\vr^{0}(s,\tilde X^{0}_{s})]^{2}ds \right]\\
 &=\frac{\epsilon\lambda_{2}}{4\lambda_{1}}\E\left[  (\partial_x \hat g(\tilde X^{0}_{T})-\partial_x \vr^{0}(0,x))^{2}\right]
 \\
 &=\frac{\epsilon\lambda_{2}}{4\lambda_{1}}\E\left[  (\partial_x \hat g(\tilde X^{0}_{T})-\E[\partial_x \hat g(\tilde X^{0}_{T})])^{2}\right]\\
 &=\frac{\epsilon\lambda_{2}}{4\lambda_{1}}{\rm Var}\left[  \partial_x \hat g(\tilde X^{0}_{T})\right]
 \end{align*}
 and the computation of the gamma $\partial_x ^{2}\vr^{0}$ is not required. Such a formulation does not seem available in general. 
 \end{remarque}
 
 \def\Xe{X^{\epsilon}}
 \def\X0{\tilde X^{0}}
 \def\sige{\sigma_{\epsilon}}
 \def\Fe{F_{\epsilon}}

 The expansion of Proposition \ref{prop: asymptotique} leads to a natural approximate hedging strategy. The result is stated in terms of the function $\Delta v$ introduced in the proof of Proposition \ref{prop: asymptotique}, see \eqref{eq: Delta vr}.
 \begin{prop}\label{prop: approximate hedging} Assume that the conditions of Proposition \ref{prop: asymptotique} hold  and that  
 \begin{itemize}
 \item[\rm (i)] $\partial^{2}_{z}\bar F_{0}\in C^{1,2}_{b}([0,T]\x \R)\cap C^{0,4}_{b}([0,T]\x \R)$,
 \item[\rm (ii)] $(t,x,z)\in [0,T]\x \R\x \R \mapsto \frac1{2\epsilon} \sigma^{2}(t,x, \epsilon z)$ is bounded and uniformly Lipschitz in its two last components, uniformly in $\epsilon \in (0,\epsilon_{0}]$.
 \end{itemize}
  Then, there exists  a constant $C>0$    such that, for each $\epsilon\in  (0,\epsilon_{0}]$  and $x\in \R$,   
 $$
 |V^{\epsilon,0,x,v^{\epsilon},\phi^{\epsilon}}_{T}- \hat g(X^{\epsilon,0,x,\phi^{\epsilon}}_{T})|\le C \epsilon^{2}
 $$
 in which 
 $$
 v^{\epsilon}:=\vr^{0}(0,x)+\epsilon \Delta v(0,x)
 $$
 and $\phi^{\epsilon}=(y^{\epsilon},b^{\epsilon},\gamma^{\epsilon})\in \Ac$ with 
 \begin{align*}
 y^{\epsilon}&=\partial_x (\vr^{0}+\epsilon \Delta v)(0,x),\\
 b^{\epsilon}&=\left[\partial_{ t}+\frac1{2\epsilon} \sigma^{2}(\cdot, \epsilon \partial_x ^{2}(\vr^{0}+\epsilon \Delta v)) \partial_x ^{2}\right]  \partial_x (\vr^{0}+\epsilon \Delta v)(\cdot, X^{\epsilon,0,x,\phi^{\epsilon}}),\\
 \gamma^{\epsilon}&= \partial_x ^{2}(\vr^{0}+\epsilon \Delta v)(\cdot, X^{\epsilon,0,x,\phi^{\epsilon}}).   
 \end{align*}
  \end{prop}
  
 \proof For ease of notations, we write $\sigma_{\epsilon}$ for $\epsilon^{-\frac12}\sigma(\cdot,\epsilon \cdot)$. We let 
 $ Y^{\epsilon}= \partial_x (\vr^{0}+\epsilon \Delta v)(\cdot, X^{\epsilon,0,x,\phi^{\epsilon}})$,  and only write $X^{\epsilon}$ for $X^{\epsilon,0,x,\phi^{\epsilon}}$ in the following. 
 Note that    \eqref{eq: tilde Xz}, \eqref{eq: Delta vr}, (i) and \eqref{eq: bar F unif elliptic} imply that $\Delta v\in C^{1,2}_{b}([0,T]\x \R)\cap C^{0,4}_{b}([0,T]\x \R)$.
 Then, the dynamics are well-defined thanks to Remark \ref{rem: vr 0 est smooth}, and  $\phi^{\epsilon}\in \Ac$. Set $F_{\epsilon}:=F(\cdot,\epsilon\;\cdot)/\epsilon$. 
By applying It\^{o}'s Lemma, using Remark \ref{rem: vr 0 est smooth}, Remark \ref{rem: edp Delta v} and  the definition of $\bar F_{\epsilon}$ together with  \eqref{eq: bar F(,0)=0}, we obtain
 \begin{align*}
 &\hat g(\Xe_{T})-v^{\epsilon}-\int_{0}^{T} Y^{\epsilon}_{t}dX^{\epsilon}_{t}-\int_{0}^{T}  \Fe(t,\Xe_{t},\gamma^{\epsilon}_{t})dt \\
 =&\vr^{0}(T,\Xe_{T})+\epsilon \Delta v(T,\Xe_{T})-\vr^{0}(0,x)-\epsilon \Delta v(0,x) -\int_{0}^{T} Y^{\epsilon}_{t}dX^{\epsilon}_{t}\\
 &-\int_{0}^{T}  \Fe(\cdot,\partial_x ^{2}(\vr^{0}+\epsilon \Delta v))(t,\Xe_{t})dt\\
 =& \int_{0}^{T} \left[\bar  \Fe(\cdot,\partial_x ^{2}(\vr^{0}+\epsilon \Delta v)) - \partial_{z}\bar F_{0}\partial_x ^{2}(\vr^{0}+\epsilon \Delta v) - \frac\epsilon 2  \partial^2_z  \bar F_0 (\partial_x ^{2}\vr^{0})^{2}\right] (t,\Xe_{t})dt.
  \end{align*}
  Recalling that \eqref{eq : regul bar F cas general} is assumed to hold for $\bar \Fe$, uniformly in $\epsilon\in (0,\epsilon_{\circ}]$,   that $\dxx \vr^{0}$ and $\dxx \Delta v$ are bounded, as well as \eqref{eq: bar F(,0)=0},   a second order  Taylor expansion implies     
  $$
  \bar  \Fe(\cdot,\partial_x ^{2}(\vr^{0}+\epsilon \Delta v)) - \partial_{z}\bar F_{0}\partial_x ^{2}(\vr^{0}+\epsilon \Delta v) - \frac\epsilon 2  \partial^2_z  \bar F_0 (\partial_x ^{2}\vr^{0})^{2}=O(\epsilon^{2}),
  $$
  in which $O(\epsilon^{2})$ is uniform on $[0,T]\x \R$.   
 \endp

 \section{Dual representation formula in the convex case} \label{sec: dual}
 In this last section, we assume that 
\begin{align}
&z\in \R\mapsto \bar F(t,x,z)\;\mbox{ is convex  and  bounded from below, }\label{eq: bar F convex}\\
& \lim_{z\to \bar \gamma(t,x)} \partial_{z}\bar F(t,x,z)=\infty\;\;\mbox{ for all $(t,x)\in [0,T]\x \R$.} \label{eq: bar F verfifie inada en l infini}
 \end{align}
  
Note that the second assumption is automatically satisfied if $\bar \gamma<\infty$, since in this case  $\lim_{z\to \bar \gamma(t,x)} \bar F(t,x,z)=\infty$. 
Both are satisfied is the model studied in  \cite{BoLoZo2}, see Remark \ref{rem: ex bar F}. \\

Whenever $\bar \gamma<\infty$, let us now use the extension $\bar F(\cdot,z):=\infty$ for $z\in [\bar \gamma,\infty)$ and define the Fenchel-Moreau transform  
$$
\bar F^{*}(\cdot,\mathfrak v):=\sup_{z\in \R} \left(\frac12\mathfrak v z-\bar F(\cdot,z)\right),\;\mathfrak v \in \R.
$$

The conditions \eqref{eq: bar F convex} and \eqref{eq: bar F verfifie inada en l infini} ensure that  $\bar F^{*}(t,x,\cdot)$ is finite on $\R_{+}$ and takes the value $+\infty$ on $\R_{-}\setminus \{0\}$. The function $\bar F$ being lower-semicontinuous on $\R_{+}$, convex and proper in its last argument, it follows that 
\begin{align}\label{eq: fenchel duality}
\bar F(\cdot,z)&=\sup_{\sr \in \R_{+}} \left(\frac12 \sr^{2} z-\bar F^{*}(\cdot,\sr^{2})\right). \\
\label{eq: optimum dans Fenchel}
\bar F^{*}(\cdot,2\partial_{z}\bar F(\cdot,z))&= \partial_{z}\bar F(\cdot,z) z-\bar F(\cdot,z),\;\mbox{ for }z<\bar \gamma.
\end{align}

 \begin{remarque}\label{rem : equiv pde} 
It follows from \eqref{eq: fenchel duality} that a function $V$ is a viscosity supersolution (resp.~subsolution) on $[0,T)\x \R$ of 
 \begin{align*}
 &\min\{-\partial_{t}\vp- \bar F(\cdot,\partial_x ^{2}\vp)\;,\; \bar \gamma - \partial_x ^{2}\vp\}=0  
 \end{align*}
 if and only if it is a viscosity supersolution (resp.~subsolution) on $[0,T)\x \R$ of 
  \begin{align}\label{eq: PDE F*}
 \inf_{\sr \in \R_{+}}\left(\bar F^{*}(\cdot,\sr^{2})  -\partial_{t}\vp-\frac12 \sr^{2} \partial_x ^{2}\vp\right) =0.
 \end{align}
 \end{remarque}
 
 This suggests, in the spirit of \cite{SonTouzZhang}, that $\vr$ admits a dual formulation in terms of an optimal control problem.  
 
 \begin{theo}\label{thm: dual formulation}  Assume   that \eqref{eq: bar F convex} and \eqref{eq: bar F verfifie inada en l infini} hold.    Let $\bS$ denote the collection   of non-negative bounded  predictable processes. 
 Then, for all $(t,x)\in [0,T)\x \R$, 
 \begin{align}\label{eq: dual formulation}
 \vr(t,x)&=\sup_{\smf \in \bS} \E\left[\hat g(X^{t,x,\smf}_{T})-\int_{t}^{T} \bar F^{*}(s,X^{t,x,\smf}_{s},\smf_{s}^{2})ds\right]\\
 &=\sup_{\smf \in \bS} \E\left[  g(X^{t,x,\smf}_{T})-\int_{t}^{T} \bar F^{*}(s,X^{t,x,\smf}_{s},\smf_{s}^{2})ds\right]\nonumber
 \end{align}
 in which 
  $$
 X^{t,x,\smf}=x+\int_{t}^{\cdot} \smf_{s}dW_{s}, \; \mbox{ $\smf \in \bS$.}
 $$  
 If moreover the conditions of Corollary \ref{cor: vr=u for Phi=g} hold, then the optimum is achieved by the Markovian control 
 $$
\hat \smf_{t,x} := (2\partial_{z}\bar F(\cdot, \partial_x ^{2}\vr)(\cdot, X^{t,x,\hat \smf_{t,x}}))^{\frac12}.
$$
 \end{theo}
 
 \begin{remarque} The model studied in \cite{BoLoZo2} corresponds to 
\begin{align*}
\bar F^{*}(t,x,s^{2})=\frac12 \frac{(s-\sigma_\circ(t,x))^{2}}{f(x)},\;\mbox{ for $s\ge 0$.}
\end{align*}
See Remark \ref{rem: ex bar F}.
The result of Theorem \ref{thm: dual formulation} above can then be formally interpreted as follows. The larger the impact function $f$, the more the optimal control can deviate from the  volatility associated to the model without market impact. When $f$ tends to $0$,   the optimal control needs to converge to  the  volatility of the impact free model $\sigma_\circ$, and one recovers the usual pricing rule at the limit. 
 \end{remarque}
 
\noindent {\bf Proof of Theorem \ref{thm: dual formulation}.} 1. We first prove the first equality in \eqref{eq: dual formulation} in the case where the conditions of   Corollary \ref{cor: vr=u for Phi=g} hold. 
Let $v$ denote the right-hand side of \eqref{eq: dual formulation}.  Recalling from Remark \ref{rem : equiv pde}, Corollary \ref{cor: vr=u for Phi=g} and Theorem \ref{theo: u C4} that $\vr$ is a smooth supersolution of \eqref{eq: PDE F*}, we deduce that 
$\vr \ge    v$  by a simple verification argument. 
Let now $\hat X$ be the solution of  
$$
  \hat  X=x+\int_{t}^{\cdot} (2\partial_{z}\bar F(\cdot, \partial_x ^{2}\vr)(s,\hat X_{s}))^{\frac12} dW_{s}.
$$  
It is well defined, recall Corollary \ref{cor: vr=u for Phi=g}, Theorem \ref{theo: u C4}, \eqref{eq: bar F unif elliptic} and \eqref{eq : regul bar F cas general}, and corresponds to $X^{t,x,\hat \smf}$ with  
$$
\hat \smf := (2\partial_{z}\bar F(\cdot, \partial_x ^{2}\vr)(\cdot,\hat X ))^{\frac12},
$$
which is bounded. Moreover, \eqref{eq: optimum dans Fenchel} implies that 
\ben
\vr(t,x) = \E\Big[\hat g(\hat  X_T) - \int_t^T\bar{F}^*(s,\hat  X_{s},{\hat \smf }_{s}^2)ds\Big],
\enn
which shows that $\vr\le   v$ since $\hat \smf$ is bounded. 

2. We now extend   the first equality in \eqref{eq: dual formulation} to the general case. Let $\{\Phi\epu,\Phi\ep\}$ be as in the proof of   Theorem \ref{theo: vr=u for Phi=hat g 0} at the end of Section \ref{sec: pde cara}, and let   $u\ep$ and $u\epu$ be the (smooth) solutions to {\rm (Eq$_{0}$)} associated to $\Phi\ep$ and $\Phi\epu$ respectively, as in Theorem \ref{theo: vr=u for Phi=hat g 0}. Then $\Phi\epu\le \hat g\le \Phi\ep$, $u\epu\le \vr \le u\ep$ and $(u\ep-u\epu, \Phi\ep-\Phi\epu)_{\eps>0}$ converges uniformly to $0$ as $\eps\to 0$. Define  $v\epu$ and $v\ep$ as $v$ but  with $\Phi\epu$ and $\Phi\ep$ in place of $\hat g$. Then, $v\epu\le v\le v\ep$ and $(v\ep-v\epu)_{\eps>0}$ converges uniformly to $0$ as $\eps\to 0$. Since, by 1., $(v\epu,v\ep)=(u\epu,u\ep)$, the required result follows.

3. It remains to prove the second equality in \eqref{eq: dual formulation}.  Define 
$$
\tilde v(t,x):=\sup_{\smf \in \bS} \E\left[  g(X^{t,x,\smf}_{T})-\int_{t}^{T} \bar F^{*}(s,X^{t,x,\smf}_{s},\smf_{s}^{2})ds\right], \;(t,x)\in [0,T)\x \R.
$$
In view of 2., we know that $\tilde v$ is bounded from above by $\vr$. Since $\bar F^{*}(\cdot,0)^{+}$ and $g^{-}$ are bounded, see \eqref{eq: bar F convex} and \eqref{eq: hyp g}, it is also bounded from below, by a constant. 
Then, it follows from \cite{bouchard2011weak} that the lower-semicontinuous enveloppe $\tilde v_{*}$ of $\tilde v$ is a viscosity supersolution of  \eqref{eq: PDE F*} such that $\tilde v_{*}(T,\cdot)\ge g$, recall \eqref{eq: hyp g}. It is in particular a supersolution of $\bar \gamma-\dxx \vp\ge 0$ on $[0,T)\x \R$, by Remark \ref{rem : equiv pde}. Then, the same arguments as in \cite[Step 3.b., proof of Theorem 3.16]{BoLoZo2} imply that $\tilde v_{*}(T,\cdot)\ge \hat g$.  By \cite{bouchard2011weak} again, we also have that 
\begin{align*}
\tilde v(t,x)&\ge \E\left[  \tilde v_{*}(T,X^{t,x,\smf}_{T})-\int_{t}^{T} \bar F^{*}(s,X^{t,x,\smf}_{s},\smf_{s}^{2})ds\right],\;\mbox{ for any $\smf \in \bS$.}
\end{align*}
Hence, 
\begin{align*}
\tilde v(t,x)&\ge \sup_{\smf \in \bS} \E\left[  \hat g(X^{t,x,\smf}_{T})-\int_{t}^{T} \bar F^{*}(s,X^{t,x,\smf}_{s},\smf_{s}^{2})ds\right].
\end{align*}

\endp
\\

We conclude this section with a result showing that any optimal control control $\hat \smf$  should be such that $\hat g(X^{t,x,\hat \smf}_{T})=g(X^{t,x,\hat \smf}_{T})$.  

\begin{prop} Let the condition of Theorem \ref{thm: dual formulation} hold and assume that $\bar F(\cdot, \kappa)$ is uniformly bounded on $[0,T]\x \R$ for some $\kappa>0$. Fix $(t,x)\in [0,T)\x \R$ and let $(\smf^{n})_{n\ge 1}$ be such that 
$$
\vr(t,x)=\lim_{n\uparrow \infty} \E\left[   g(X^{t,x,\smf^{n}}_{T})-\int_{t}^{T} \bar F^{*}(s,X^{t,x,\smf^{n}}_{s},(\smf^{n}_{s})^{2})ds\right].
$$
Then, $(X^{t,x,\smf^{n}}_{T})_{n\ge 1}$ is tight, and any limiting law $\nu$ associated to a subsequence satisfies 
$
\nu(\hat g>g)=0.
$
\end{prop}

\proof We only write $X^{n}$ for $X^{t,x,\smf^{n}}$ and let $$J_{n}:= \E\left[   g(X^{{n}}_{T})-\int_{t}^{T} \bar F^{*}(s,X^{^{n}}_{s},(\smf^{n}_{s})^{2})ds\right],$$ $n\ge 1$. 
Then, \eqref{eq: hyp g} and \eqref{eq: bar F convex} imply that one can find $C>0$ such that 
\begin{align*}
-C&\le \E[C+\frac{\kappa}4 |X^{n}_{T}|^{2}-\int_{t}^{T} \frac{\kappa}2 (\smf^{n}_{s})^{2}ds + T\sup \bar F(\cdot, \kappa)]\\
&\le \E[C-\int_{t}^{T} \frac{\kappa}4 (\smf^{n}_{s})^{2}ds + T\sup \bar F(\cdot, \kappa)].
\end{align*}
Hence, $\sup_{n\ge 1}\E[\int_{t}^{T}   (\smf^{n}_{s})^{2}ds]<\infty$. Let $\nu_{n}$ be the law associated to $X^{n}_{T}$. The above shows that $(\nu_{n})_{n\ge 1}$ is tight. Let us consider a subsequence $(\nu_{n_{k}})_{k\ge 1}$ that converges to   some law $\nu$.  If $\nu(\hat g>g)>0$, then one can find $\delta>0$ such that 
$\E[\hat g(X^{n_{k}}_{T})]\ge \E[g(X^{n_{k}}_{T})]+\delta$ for all $k\ge 1$ large enough,  which would imply that 
$$
\lim_{k\to  \infty} \E\left[   \hat g(X^{{n_{k}}}_{T})-\int_{t}^{T} \bar F^{*}(s,X^{{n_{k}}}_{s},(\smf^{n_{k}}_{s})^{2})ds\right]\ge \lim_{k\to  \infty} J_{n_{k}}+\delta,
$$
a contradiction to Theorem \ref{thm: dual formulation}.\endp
\\ 

\bibliography{biblio-greg-1}
\end{document}